\documentclass[11pt,a4paper]{article}
\usepackage[latin1]{inputenc}
\usepackage[american]{babel}      
\usepackage[T1]{fontenc}
\usepackage{amsmath,amssymb,amsthm}
\usepackage{graphicx}
\usepackage{xcolor}   
\usepackage{geometry}

\title{{\bf Saturation theorems for neural network operators by solving elliptic and hyperbolic PDEs with analytical and semi-analytical inverse problems
}}
         
\author{ {\bf Danilo Costarelli \thanks{Paper submitted to arXiv: arXiv:2505.07526v1}} \\  
Department of Mathematics and Computer Science \\
            University of Perugia\\
        1, Via Vanvitelli, 06123 Perugia, Italy    \\  
{\small {\tt danilo.costarelli@unipg.it}} }

\date{}

\newcommand{\mau}{\geq}
\newcommand{\miu}{\leq}

\newcommand{\N}{\mathbb{N}}
\newcommand{\R}{\mathbb{R}}
\newcommand{\Z}{\mathbb{Z}}
\newcommand{\disp}{\displaystyle}
\newcommand{\be}{\begin{equation}}
\newcommand{\ee}{\end{equation}}
\newcommand{\phis}{\phi_{\sigma}}
\newcommand{\psis}{\Psi_{\sigma}}

\newcommand{\xx}{{\tt x}}
\newcommand{\yy}{{\tt y}}

\newtheorem{definition}{Definition}[section]
\newtheorem{remark}[definition]{Remark}
\newtheorem{theorem}[definition]{Theorem}
\newtheorem{lemma}[definition]{Lemma}

\newtheorem{example}[definition]{Example}
\newtheorem{corollary}[definition]{Corollary}

\begin{document}

\maketitle 

\begin{abstract}
This paper addresses inverse problems (in a broad sense) for two classes of multivariate neural network (NN) operators, with particular emphasis on saturation results, and both analytical and semi-analytical inverse theorems. One of the key aspects in addressing these issues is solving certain elliptic and hyperbolic partial differential equations (PDEs), as well as suitable asymptotic formulas for the NN operators based on sufficiently smooth functions; the connection between these two topics lies in the application of the so-called generalized parabola technique by Ditzian. 
The saturation theorems yield characterizations of the corresponding saturation classes; these are respectively related to harmonic functions and to the solution of a certain transport equation. Analytical inverse theorems, on the other hand, are related to sub-harmonic functions as well as to functions in the Sobolev space $W^2_\infty$. Finally, we also address a semi-analytical inverse problem. The latter serves as the starting point for deriving a retrieval procedure that may be useful in real world applications.

\vskip0.3cm
\noindent
  {\footnotesize AMS 2010 Mathematics Subject Classification: 35A02, 31B20, 32W50, 35A09, 35B65, 41A25}
\vskip0.1cm
\noindent
  {\footnotesize Key words and phrases: inverse problems; saturation classes, elliptic PDEs; hyperbolic PDEs; harmonic functions; transport equation; retrieval} 
\end{abstract}

\section{Introduction} \label{sec1}

 Over the past ten years, the theory of the so-called neural network (NN) operators activated by sigmoidal functions $\sigma$ has been extensively studied, mainly with regard to their convergence properties and their order of approximation. A wide selection of recent papers on this topic includes the following \cite{AN-2017,BAKU-2021,QIYU2022,KAD-2023,KACO-2023,Acar-2024,BAX24,DEM-2024,SINGH,MABA-2025,YUCA2025}.

Delicate issues (especially when one deals with multivariate data) concerning inverse problems have not previously been studied; in particular we refer to saturation results, and both analytical and semi-analytical inverse theorems.

This is the main motivation behind the present paper, which is entirely devoted to {\em inverse problems} for the two main studied families of NN operators, that are the classical NN operators (denoted in what follows by $F_n$) introduced in \cite{COSP3} and the Kantorovich NN operators (denoted by $K_n$) given in \cite{COSP4}. Concerning the main motivations that led to the introduction of the so-called NN operators, we directly refer to \cite{CAEU1992,ANAST1,CAOCHEN2009,COSP2,YU25,MOW} where also the relation with artificial intelligence and machine learning have been discussed.

By a {\em saturation result} (\cite{DIMA1976}), we mean establishing the best possible non-trivial order of approximation that can be achieved by a given linear or nonlinear net of operators; at the same time, if a function is approximated at an order higher than the saturation order, we expect that they belong to a certain special class which characterizes the considered approximation process. 

Analytical inverse theorems, on the other hand, aim to determine the regularity of a given function, knowing the order of approximation with which it can be approximated; this is exactly the converse problem of what is generally known as direct approximation.

One of the key aspects in addressing the above issues is the solution of certain elliptic and hyperbolic partial differential equations (PDEs), as well as suitable asymptotic formulas for NN operators based on sufficiently smooth functions; the link between the two mentioned topics is the application of the so-called generalized parabola technique by Ditzian (\cite{Ditzian1983}). 
The above non-trivial relation between PDEs and inverse problems has been partially analyzed in \cite{CARDENAS} in the general context of simultaneous approximation by positive linear operators, mainly in the case of relations with elliptic problems, in which the authors strongly use the maximum principle for elliptic PDEs (\cite{GILB}). To the best of our knowledge, the relations between saturation results and hyperbolic PDEs have never been considered until now.

In addition to the analytical results, we also address the problem of the semi-analytical inversion of operators. This is in general a nonlinear problem which is strictly related to real-world applications; in order to be solved it needs both analytical and numerical methods. Indeed, in certain real-world applications (for instance, in remote sensing data analysis, see, e.g., \cite{review1}) it may happen that, when certain data are detected, one may be interested in the variables influencing these measurements. This can be important since certain functional dependencies are theoretically known but cannot be explicitly determined.

Based on the above considerations, it could be useful to have a procedure that, starting from the measured data, allows us to estimate the value of certain target variables; this type of inverse problem is usually referred to as a {\em retrieval problem} (\cite{review2}).

Modelling by NN operators provides (approximate) functional expressions linking a number of variables to scalar values, even when the explicit functional dependencies are unknown. Hence, their semi-analytical inversion can be useful to give a possible solution of the retrieval problem.
\vskip0.2cm

Below we provide a more detailed description of all the results established in this paper and the considered techniques of proof. 
\vskip0.2cm

  After a brief review of the foundations of the theory of NN operators (see Section \ref{sec2}), in Section \ref{sec3} we immediately considered analytical inverse problems. For both the multivariate NN operators $F_n$ and $K_n$ we establish an inverse result which asserts that, if a function  can be approximated by an error that with respect to the uniform norm decays as the power $n^{-\nu}$, $0<\nu<1$, then $f$ belongs to the Lipschitz (or Holderian) space $Lip(\nu)$. Thanks to this result, also exploiting the direct approximation theorems established in \cite{COCOKA1}, we achieve a characterization of $Lip(\nu)$ in terms of the order of approximation of the NN operators. 

  From the latter result, the natural question concerning what is the saturation order for the NN operators arises. By a simple example it is possible to see that the best possible {\em global} order of approximation is ${\cal O}(n^{-1})$, as $n \to +\infty$, and it can not be improved in general. At the same time, also a two-sided estimate for the constant of best approximation can be achieved, i.e., an upper and lower bound for the constants used for estimating the decay of the approximation error can be derived. The latter issue is important since the {\em optimal constants} of approximation are fundamental to evaluate the performance of an approximation method.

In Section \ref{sec4}, local asymptotic formulas have been proved for $F_n$ and $K_n$. In order to do this, using methods of Fourier Analysis, certain sums involving the kernel functions defining the considered operators (i.e., the so-called {\em truncated algebraic moments}) have been estimated. In particular, we proved that under suitable non-restrictive assumptions, any truncated algebraic moment (t.a.m) of odd degree vanishes pointwise, while those with even order are all constants and positive. In particular, the second order t.a.m. can be computed exploiting the second derivative of the Fourier transform of the kernel function evaluated at zero: this is a consequence of the application of the Poisson summation formula to a specific 1-periodic function defined by a certain series. From the mentioned asymptotic analysis, we deduce that by the operators $F_n$ it is possible to establish local approximations that are of the order ${\cal O}(n^{-2})$, as $n \to +\infty$. While, concerning the operators $K_n$ local approximations have lower order than those achieved by $F_n$. The latter fact is not surprising since the operators $F_n$ are defined exploiting the pointwise values of the reconstructed function $f$, while the operators $K_n$ are defined by means of certain averages of the function $f$; thus, the first are "by definition" more accurate than the second one, that seems to be more effective in case of not necessarily continuous functions and to contrast the possible perturbation effects on the sampled $f$.  

The results achieved in Section \ref{sec4} are fundamental to achieve a characterization of the so-called {\em saturation classes}, i.e., the non-trivial sub-classes of the space of continuous functions that are characterized by the best possible order of approximation.

In detail, in Section \ref{sec5} we considered the latter problem for the operators $F_n$. By a suitable revisiting of the generalized parabola method introduced by Ditzian (\cite{Ditzian1983}) we prove that the (local) saturation classes are provided by harmonic functions; in the proof of such a result a Dirichlet problem for the Laplace equation is solved. By a similar strategy we also derive the inverse result for the case corresponding to the local saturation order; in this case we find a certain modification of $f$ that turns out to be a sub-harmonic function. In the special case of the univariate functions, it turns out that the saturation class corresponds to the Sobolev space $W^2_\infty$.

The results just mentioned are only {\em local}; in the one-dimensional case $d=1$ also the global saturation order has been established, and we proved that the only functions that can be approximated by a higher order are the constant ones.

  Here we also present concrete examples of NN operators, for instance the one based on the logistic function and on the sigmoidal function generated by the hyperbolic tangent, as well as suitable sigmoidal functions generated by an integral function involving the so-called Jackson-type kernels.

Concerning the operators $K_n$ instead (see Section \ref{sec6}), we prove a local theorem which is instead connected to the solution of a certain hyperbolic PDE. More precisely, we are referring to a transport equation with constant coefficients (or the scalar homogeneous divergence equation). We stress again that to the best of our knowledge, this connection between saturation theorems and hyperbolic PDEs has not been studied.


Finally, in Section \ref{sec8}, the above mentioned semi-analytical inverse problem has been faced. In the literature, we can find different approaches for the inversion of operators, which use both methods of Functional and Numerical Analysis. For instance, we can mention the approach of \cite{HOWLET} based on the expansion of the considered operators using the McLaurin series. Here, we exploit the expression of a certain sigmoidal function known as the {\em ramp function} to combine analytical and approximate methods, justifying the used name of semi-analytical approach. To show the concrete usefulness of the proposed retrieval procedure, we also provided two concrete examples of retrieval for both functions of one- and two- variables.
%


\section{Some recall on the multivariate neural network operators} \label{sec2}

A function $\sigma: \R \to \R$ is said to be a {\em sigmoidal function} if $\lim_{x \to -\infty} \sigma(x) = 0$ and $\lim_{x \to +\infty} \sigma(x) = 1$ (\cite{CYB}). 

Let now $\sigma$ be a given non-decreasing sigmoidal function, such that:
\begin{description}
	\item[$(\Sigma 1)$] $g_{\sigma}(x) := \sigma(x) - 1/2$ is an odd function;
	\item[$(\Sigma 2)$] $\sigma \in C^2(\R)$ is concave for $x \geq 0$;
	\item[$(\Sigma 3)$] $\sigma(x) = \mathcal{O}(|x|^{-1-\alpha})$ as $x \to -\infty$, for some $\alpha > 0$,
\end{description}
are satisfied. Under conditions $(\Sigma 1)$ and $(\Sigma 3)$ one has:
\be \label{ordine-simmetrico}
1-\sigma(x)\ =\ \sigma(-x)\ =\ \mathcal{O}(|x|^{-1-\alpha}), \quad\ as \quad x \to +\infty.
\ee
Hence, we define
\begin{equation*}
   \phi_{\sigma}(x) := \frac{1}{2}[\sigma(x+1) - \sigma(x-1)], \quad  x \in \R,
\end{equation*}
a density function generated by $\sigma$ (for more details, see \cite{COSP2}).
We now recall the definition of the multivariate density function generated by $\sigma$:
\begin{equation*}
  \Psi_{\sigma}({\tt{x}}) := \phi_{\sigma}(x_1) \cdot \phi_{\sigma}(x_2) \cdot \dots \cdot \phi_{\sigma}(x_d), 
  \quad  {\tt{x}} := (x_1, ..., x_d) \in \R^d.	
\end{equation*} 
For some basic properties of $\Psi_{\sigma}$, see, e.g., \cite{COSP3}.
\vskip0.2cm

From now on, when we use the symbol $\sigma$ we always refer to a sigmoidal function as above, satisfying conditions $(\Sigma 1)$, $(\Sigma 2)$ and $(\Sigma 3)$. In what follows we denote by $I:=\prod_{i=1}^{d} [a_i,b_i]$ the $d$-dimensional set on $\R^d$, with $a_i,b_i\in\R$, $i=1,\cdots,d$, and by $C(I)$ the space of all continuous real-valued functions defined on $I$, equipped with the sup-norm $\| \cdot \|_{\infty}$. Further, let introduce the following set
\begin{equation*}
\mathcal{J}_n:=\{{\tt{k}}\in\Z^d \;:\; \lceil {na_i \rceil } \leq k_i \leq \lfloor {nb_i \rfloor},\; i=1,\cdots,d \}, \quad n\in\N,
\end{equation*}
where the symbols $\lceil \cdot \rceil$ and $\lfloor \cdot \rfloor$ denote the ceiling and the floor functions, respectively.
In \cite{COSP3}, the following property has been established, for ${\tt{x}}\in I$ and $n\in\N^+$:
\begin{equation} \label{denominator1}
   m_{0,0}^n(\Psi_\sigma,n{\tt{x}}):= \sum_{{\tt{k}}\in \mathcal{J}_n}\Psi_\sigma(n{\tt{x}}-{\tt{k}})\geq \left[\phis(1)\right]^d>0,
\end{equation}
where $\phis(1)>0$ by the properties of $\sigma$. If we also define the set:
\begin{equation*}
\mathcal{K}_n:=\{{\tt{k}}\in\Z^d \;:\; \lceil {na_i \rceil } \leq k_i \leq \lfloor {nb_i \rfloor}-1,\; i=1,\cdots,d \}, \quad n\in\N,
\end{equation*}
and we assume that $\sigma(1)<1$ (see Remark 3.3 of \cite{CP1}), we can also have:
\begin{equation} \label{denominator2}
   \widetilde{m}_{0,0}^n(\Psi_\sigma,n{\tt{x}}):= \sum_{{\tt{k}}\in \mathcal{K}_n}\Psi_\sigma(n{\tt{x}}-{\tt{k}})\geq \left[\phis(2)\right]^d>0.
\end{equation}
We now recall the definitions of some very useful tools.

  We define the \textit{discrete absolute moment of $\Psi_{\sigma}$ of order $\beta\geq0$} as
\begin{equation}
M_{\beta}(\Psi_{\sigma}):=\sup_{{\tt{x}}\in\R^d}\sum_{{\tt{k}}\in \Z^d} \Psi_\sigma({\tt{x}}-{\tt{k}})\left\|{\tt{k}}-{\tt{x}}\right\|^{\beta}_2,
\end{equation}
where $\left\|\cdot\right\|_2$ denotes the usual Euclidean norm on $\R^d$.
It is well-known that (see \cite{COSP3}) 
\begin{equation}\label{momento}
M_0(\Psi_\sigma)=1, \quad and \quad   M_{\beta}(\Psi_{\sigma})<+\infty, \qquad 0\leq \beta<\alpha, \quad \quad d \mau 1,
\end{equation}
where for $d=1$ the above properties hold for $\phis$ in place of $\psis$. 
\vskip0.2cm

Now, we recall respectively the definitions of the classical and the Kantorovich neural network (NN) operators, that will be considered in this paper.
\vskip0.2cm 

  Let $f: I \to \R$ be a given function, and $n \in \N^+$. The 
NN {\em operators} $F_n(f,\xx)$, activated 
by the sigmoidal function $\sigma$ (acting on $f$) are defined by
\begin{equation*}
	 F_n(f,\xx) := \frac{\disp \sum_{{\tt{k}}\in \mathcal{J}_n}f\left({\tt{k}}/n \right) \Psi_{\sigma}(n{\tt{x}}-{\tt{k}})}
      {\disp \sum_{{\tt{k}}\in \mathcal{J}_n}
           \Psi_{\sigma}(n{\tt{x}}-{\tt{k}})}, 
\end{equation*}
for every $\xx \in I$ and ${\tt{k}}/n := (k_1/n, ..., k_d/n)$. Furthermore, the Kantorovich NN operators activated by $\sigma$ are
\begin{equation*}
	 K_n(f,\xx) := \frac{\disp \sum_{{\tt{k}}\in \mathcal{K}_n}\left[  n^d \int_{{\tt{k}}/n}^{({\tt{k}+1})/n} f({\tt u})\, d{\tt u}   \right]\, \Psi_{\sigma}(n{\tt{x}}-{\tt{k}})}
      {\disp \sum_{{\tt{k}}\in \mathcal{K}_n}
           \Psi_{\sigma}(n{\tt{x}}-{\tt{k}})}, 
\end{equation*} 
where $f$ is locally integrable, and by the integrals:
$$
\int_{{\tt{k}}/n}^{({\tt{k}+1})/n} f({\tt u})\, d{\tt u},
$$
we briefly denote the multiple integrals on the multivariate rectangles:
$$
\left[{k_1 \over n}, {k_1+1 \over n}  \right] \times ... \times \left[{k_d \over n}, {k_d+1 \over n}  \right]\ \subset\ I.
$$
Sometimes, to denote the above operators we also use the equivalent (and simpler) notation $F_nf$ and $K_nf$.

  Both operators are well-defined thanks to (\ref{denominator1}) and (\ref{denominator2}) respectively, provided, for instance, that the function $f$ is bounded. Indeed, it is easy to see that 
$$
\left|F_n(f,{\tt{x}})\right|\leq \|f\|_\infty, \quad \left|K_n(f,{\tt{x}})\right|\leq \|f\|_\infty, \quad {\tt x} \in I, \quad n \in \N.
$$

For the study of the convergence properties of the above operators we refer to \cite{COSP3} and \cite{COSP4}, while the order of approximation has been studied in \cite{COCOKA1}.

To prove the results of next sections, it can be useful to introduce the following notations:
$$
{\cal J}_n^{[i]}\ :=\ \{{\tt{k}}\in\Z^{d-1} \;:\; \lceil {na_j \rceil } \leq k_j \leq \lfloor {nb_j \rfloor},\; j=1,\cdots,d,\, j\neq i \}, \quad n\in\N,
$$
$$
{\cal I}_n^{[i]}\ :=\ \{k_i\in\Z \;:\; \lceil {na_i \rceil } \leq k_i \leq \lfloor {nb_i \rfloor}\}, \quad n\in\N,
$$
$$
{\tt k}_{[i]}\, :=\, (k_1, .., k_{i-1},k_{i+1}, ..., k_d)\ \in \Z^{d-1}, 
$$
$$
{\tt x}_{[i]}\, :=\, (x_1, .., x_{i-1},x_{i+1}, ..., x_d)\ \in \R^{d-1}, 
$$
$$
\Psi^{[i]}_\sigma({\tt x}_{[i]})\ :=\ \phis(x_1) \dots \phis(x_{i-1})\, \phis(x_{i+1}) \dots \phis(x_d), \quad {\tt x}_{[i]} \in \R^{d-1},
$$
$i=1, ..., d$. 
\begin{remark} \rm
Note that, the (linear) operators $F_n$ and $K_n$ differ from the classical (nonlinear) deep neural network methods (see, e.g., \cite{FEZH,AVJU}) since they do not require any training and their complete analytical definition for the constructive approximation of functions is available. We also recall that $F_n$ and $K_n$ are called "neural network" operators for a main reason: this term has been coined by Cardaliaguet and Euvrard in 1992 (\cite{CAEU1992}), when the most commonly studied NNs were the feed-forward with one hidden layer. Actually, both $F_n$ and $K_n$ are feed-forward NNs activated by sigmoidal functions (except for the presence of the denominator) and although neural networks have changed considerably in recent years, the name of $F_n$ and $K_n$ has been maintained.
\end{remark}

\section{Analytical inverse results} \label{sec3}

In order to establish the main results of this section, we first recall the notion of the modulus of continuity of a given function $f:I \to \R$, as follows:
$$
\omega(f,\delta)\ :=\ \sup \left\{ |f(\xx)-f(\yy)|,\  {\tt x},\, {\tt y} \in I,\ with\ \|{\tt x}-{\tt y}\|_2 \miu \delta \right\}, \quad \delta >0.
$$
Using the above tool we can define the Lipschitz classes as:
$$
Lip(\nu)\ :=\ \{ f \in C(I):\ \omega(f,\delta)={\cal O}(\delta^\nu),\ as\ \delta \to 0^+ \}, \quad 0< \nu \miu 1.
$$
For the sake of simplicity, from now on when we refer to the symbol $\alpha$, we always refer to the parameter of assumption $(\Sigma 3)$ introduced in Section \ref{sec2}.
We can prove the following.
\begin{theorem} \label{th-inverso-NN}
Let $\alpha>1$ and assume in addition that $M_1(\phis')$ is finite. Let now $f \in C(I)$ such that:
\be \label{ipotesi-th-inverso-1}
\|   F_n(f,\cdot)-f(\cdot)\|_\infty\, =\, {\cal O}(n^{-\nu}), \quad as \quad n \to +\infty,
\ee
with $0<\nu<1$. Then $f \in Lip(\nu)$.
\end{theorem}
\begin{proof}
Let $\delta>0$ be fixed. We can write what follows:
$$
\omega(f,\delta)\ \miu\ \omega(f(\cdot)-F_n(f, \cdot),\delta)\, +\,  \omega(F_n(f, \cdot),\delta)
$$
\be \label{disuguaglianza}
\miu\ 2\, \|   F_n(f,\cdot)-f(\cdot)\|_\infty\ +\,  \omega(F_n(f, \cdot),\delta)\, \miu\ M\, n^{-\nu}\, +\, \omega(F_n(f, \cdot),\delta),
\ee
where $M>0$ is a suitable constant arising from (\ref{ipotesi-th-inverso-1}), and $n$ is sufficiently large. In order to estimate the term $\omega(F_n(f, \cdot),\delta)$ we proceed as follows. By assumption $(\Sigma 2)$ we have that the operators $F_n(f, \cdot)$ are sufficiently regular, hence the following well-known first order Taylor formula with integral remainder holds:
$$
F_n(f, {\tt u})\, -\, F_n(f,{\tt x})\, =\, \sum_{i=1}^d (u_i-x_i)\int_0^1 {\partial\, F_n f \over \partial\, x_i}({\tt x}+t({\tt u}-{\tt x}))\, dt, \quad {\tt x},\, {\tt u} \in I.
$$
Moreover, if we denote by ${\bf 1}$ the unitary constant function, we can observe that:
$$
F_n({\bf 1}, {\tt u})\ =\ 1, \quad {\tt u} \in I,
$$
hence each of its partial derivative turns out to be identically null. From this property we obtain:
$$
\hskip-9.5cm F_n(f, {\tt u})\, -\, F_n(f,\tt x)\, 
$$
$$
=\, \sum_{i=1}^d (u_i-x_i)\int_0^1 \left[ {\partial\, F_n f \over \partial\, x_i}({\tt x}+t({\tt u}-{\tt x}))\, -\,  f({\tt x}+t({\tt u}-{\tt x})){\partial\, F_n {\bf 1} \over \partial\, x_i}({\tt x}+t({\tt u}-{\tt x}))  \right]\, dt.
$$
Now, for every fixed $i=1,...,d$, computing the first order partial derivatives of $F_n f$ and $F_n {\bf 1}$, using (\ref{denominator1}) and (\ref{momento}), we can write what follows:
$$
\left|{\partial\, F_n f \over \partial\, x_i}({\tt x}+t({\tt u}-{\tt x}))\, -\,  f({\tt x}+t({\tt u}-{\tt x})){\partial\, F_n {\bf 1} \over \partial\, x_i}({\tt x}+t({\tt u}-{\tt x}))\right|
$$
$$
\miu\ \phis(1)^{-2d}\, n\,  \sum_{{\tt k} \in {\cal J}_n} \left| f({\tt k}/n)-f({\tt x}+t({\tt u}-{\tt x}))\right| \cdot  \Psi_\sigma^{[i]}\left(  n\xx_{[i]} + nt({\tt u}_{[i]}-\xx_{[i]}) -{\tt k}_{[i]} \right)\, \times
$$
$$
\left|\phis'(nx_i + nt(u_i-x_i) -k_i)\right| \, \cdot \left[ \sum_{{\tt k} \in {\cal J}_n} \psis(n{\tt x} + n t({\tt u}-{\tt x})-{\tt k})   \right]            
$$
$$
+\, \phis(1)^{-2d} \, n\, \left[\sum_{{\tt k} \in {\cal J}_n} \left| f({\tt k}/n)-f({\tt x}+t({\tt u}-{\tt x}))\right|\, \psis(n{\tt x}+n t({\tt u}-{\tt x})-{\tt k})\right]\, \times
$$
$$
\left\{   \sum_{{\tt k} \in {\cal J}_n} \Psi_\sigma^{[i]}\left(  n\xx_{[i]} + nt({\tt u}_{[i]}-\xx_{[i]}) -{\tt k}_{[i]} \right) \, |\phis'(nx_i+nt(u_i-x_i)-k_i)| \right\}
$$
$$
\miu\, \phis(1)^{-2d}\, n\,  \sum_{{\tt k} \in {\cal J}_n} \left| f({\tt k}/n)-f({\tt x}+t({\tt u}-{\tt x}))\right|\, \Psi_\sigma^{[i]}\left(  n\xx_{[i]} + nt({\tt u}_{[i]}-\xx_{[i]}) -{\tt k}_{[i]} \right)\, \times 
$$
$$
|\phis'(nx_i+nt(u_i-x_i)-k_i)| 
$$
$$
+\, \phis(1)^{-2d}\, n\, \left[\sum_{{\tt k} \in {\cal J}_n} \left| f({\tt k}/n)-f({\tt x}+t({\tt u}-{\tt x}))\right|\, \psis(n{\tt x}+nt({\tt u}-{\tt x})-{\tt k})\right]\, \times
$$
$$
\sum_{{\tt k} \in {\cal J}_n} \Psi_\sigma^{[i]}\left(  n\xx_{[i]} + nt({\tt u}_{[i]}-\xx_{[i]}) -{\tt k}_{[i]} \right)\, |\phis'(nx_i+nt(u_i-x_i)-k_i)|.
$$
Recalling the following well-known property of the modulus of continuity:
$$
\omega(f, \lambda\, \delta)\ \miu\ (1+\lambda)\, \omega(f, \delta), \quad \delta,\, \lambda >0
$$
we obtain:
$$
\left|{\partial\, F_n f \over \partial\, x_i}({\tt x}+t({\tt u}-{\tt x}))\, -\,  f({\tt x}+t({\tt u}-{\tt x})){\partial\, F_n {\bf 1} \over \partial\, x_i}({\tt x}+t({\tt u}-{\tt x}))\right|
$$
$$
\miu \phis(1)^{-2d}\, n\, \omega(f, 1/n)\,  \sum_{{\tt k} \in {\cal J}_n} \left[n\, \| {\tt k}/ n\, - {\tt x}\,+t({\tt x}-{\tt u})\|_2 + 1\right]\, \Psi_\sigma^{[i]}\left(  n\xx_{[i]} + nt({\tt u}_{[i]}-\xx_{[i]}) -{\tt k}_{[i]} \right)\, \times 
$$
$$
 |\phis'(nx_i+nt(u_i-x_i)-k_i)| 
$$
$$
+\, \phis(1)^{-2d}\, n\, \omega(f, 1/n) \left[\sum_{{\tt k} \in {\cal J}_n}  \left[n\, \| {\tt k}/ n\, - {\tt x}\,+t({\tt x}-{\tt u})\|_2 + 1\right]\,     \, \psis(n{\tt x} +nt({\tt x}-{\tt u})-{\tt k})\right]\, \times
$$
$$
\sum_{{\tt k}_{[i]} \in {\cal J}^{[i]}_n} \Psi_\sigma^{[i]}\left(  n\xx_{[i]} + nt({\tt u}_{[i]}-\xx_{[i]}) -{\tt k}_{[i]} \right)\, \sum_{k_i \in {\cal I}^{[i]}_n}|\phis'(nx_i+nt(u_i-x_i)-k_i)|.
$$
Now, denoting by:
$$
\| {\tt x}\|_m\ :=\ \max\{ |x_i|:\ i=1,...,d \}, \quad {\tt x} \in \R^d,
$$
we know $\| {\tt x}\|_2 \miu \sqrt{d}\, \| {\tt x}\|_m$, thus:
$$
\sum_{{\tt k} \in {\cal J}_n} \left[n\, \| {\tt k}/ n\, - {\tt x}\,+t({\tt x}-{\tt u})\|_2 + 1\right]\,  \Psi_\sigma^{[i]}\left(  n\xx_{[i]} + nt({\tt u}_{[i]}-\xx_{[i]}) -{\tt k}_{[i]} \right)\, |\phis'(nx_i+nt(u_i-x_i)-k_i)|
$$
$$
\miu\ \sqrt{d}\, \sum_{{\tt k} \in {\cal J}_n}\| n{\tt x}\,+nt({\tt u}-{\tt x}) - {\tt k}\|_m  \Psi_\sigma^{[i]}\left(  n\xx_{[i]} + nt({\tt u}_{[i]}-\xx_{[i]}) -{\tt k}_{[i]} \right)\, |\phis'(nx_i+nt(u_i-x_i)-k_i)|
$$
$$
+\, \sum_{k_i \in {\cal I}^{[i]}_n}|\phis'(nx_i+nt(u_i-x_i)-k_i)|
$$
$$
\leq \sqrt{d}\, \sum_{\nu = 1}^d \sum_{{\tt k} \in {\cal J}_n} |nx_\nu\,+nt(u_\nu-x_\nu) - k_\nu|  \Psi_\sigma^{[i]}\left(  n\xx_{[i]} + nt({\tt u}_{[i]}-\xx_{[i]}) -{\tt k}_{[i]} \right)\, \times 
$$
$$
|\phis'(nx_i+nt(u_i-x_i)-k_i)|\, +\, M_0(\phis') 
$$
$$
\miu\ \sqrt{d}\, M_0(\phis)^{d-2} \left[  M_1(\phis) M_0(\phis') + M_0(\phis) M_1(\phis')  \right]\, +\ M_0(\phis')
$$
$$
=\, \sqrt{d}\, \left[  M_1(\phis) M_0(\phis') + M_1(\phis')  \right]\, +\ M_0(\phis')\ =: K_1.
$$
The constant $K_1$ is independent on $n$ and it is finite since $\alpha>1$, $M_1(\phis')<+\infty$, and consequently also $M_0(\phis')<+\infty$.
Similarly:
$$
\left[\sum_{{\tt k} \in {\cal J}_n}  \left[n\, \| {\tt k}/ n\, - {\tt x}\,+t({\tt x}-{\tt u})\|_2 + 1\right]\,     \, \psis(n{\tt x}+nt({\tt x}-{\tt u})-{\tt k})\right]\, \times
$$
$$
\sum_{{\tt k}_{[i]} \in {\cal J}^{[i]}_n} \Psi_\sigma^{[i]}\left(  n\xx_{[i]} + nt({\tt u}_{[i]}-\xx_{[i]}) -{\tt k}_{[i]} \right)\,\sum_{k_i \in {\cal I}^{[i]}_n}|\phis'(nx_i+nt(u_i-x_i)-k_i)| 
$$
$$
\miu\ \left[ 1 + M_1(\psis)  \right]\, M_0(\phis')\ =:\ K_2\ <\ +\infty.
$$
Using the above constants $K_1$ and $K_2$ we finally obtain:
$$
\hskip-9.5cm |F_n(f, {\tt u})\, -\, F_n(f,\tt x)|\, 
$$
$$
\miu\, \sum_{i=1}^d |u_i-x_i|\int_0^1 \phis(1)^{-2d}\, n\,  \omega(f, 1/n)\, \left( K_1 + K_2\right) dt\ \miu\ K\, n\,  \omega(f, 1/n)\, \|{\tt u} - {\tt x}\|_2,
$$
where $K>0$ is a suitable constant. Passing to the supremum with respect to $\|{\tt u} - {\tt x}\|_2 \miu \delta$, with ${\tt u}$, ${\tt x} \in I$, and using (\ref{disuguaglianza}), we immediately get:
\be
\omega(f, \delta)\ \miu\ M\, n^{-\nu}\, +\, K\, n\,  \omega(f, 1/n)\, \delta.
\ee
Now, proceeding as in the last part of the proof of Theorem 7 of \cite{COCO1} we get the thesis.
\end{proof}
Using the same strategy, as above we can immediately state also the following result for the Kantorovich NN operators.
\begin{theorem} \label{th-inverso-KNN}
Let $\alpha>1$ and assume that $M_1(\phis')$ is finite. Let now $f \in C(I)$ such that:
$$
\|   K_n(f,\cdot)-f(\cdot)\|_\infty\, =\, {\cal O}(n^{-\nu}), \quad as \quad n \to +\infty,
$$
with $0<\nu<1$. Then $f \in Lip(\nu)$.
\end{theorem}
We omit the proof of Theorem \ref{th-inverso-KNN} since it is completely analogous to the one of Theorem \ref{th-inverso-NN}.

\vskip0.2cm

Recalling Theorem 12 and Remark 14 of \cite{COCOKA1}, we know that for sigmoidal functions with $\alpha>1$, and $f \in C(I)$ one has:
\be \label{order-NN}
\|F_n(f,\cdot)-f(\cdot)\|_\infty\ \miu\ C_1\, \omega(f, 1/n),
\ee
and 
\be \label{order-KNN}
\|K_n(f,\cdot)-f(\cdot)\|_\infty\ \miu\ C_2\, \omega(f, 1/n),
\ee
for a suitable constants $C_1>0$, $C_2>0$ and for sufficiently large $n$. Hence, combining these estimates together with Theorem \ref{th-inverso-NN} and Theorem \ref{th-inverso-KNN} we can immediately deduce the following characterization.
\begin{theorem} \label{th-carachterization}
Let $\alpha>1$; we have for $0<\nu<1$:
$$
f \in Lip(\nu) \quad \mbox{if and only if} \quad \|   F_n(f,\cdot)-f(\cdot)\|_\infty\, =\, {\cal O}(n^{-\nu}), \quad as \quad n \to +\infty,
$$
and similarly
$$
f \in Lip(\nu) \quad \mbox{if and only if} \quad \|   K_n(f,\cdot)-f(\cdot)\|_\infty\, =\, {\cal O}(n^{-\nu}), \quad as \quad n \to +\infty.
$$
\end{theorem}

Based on the above results, it is clear that for the case $\nu=1$ there are some open questions. For instance, we can ask if the {\em global} estimates in (\ref{order-NN}) and (\ref{order-KNN}) can be improved in general. On this matter we can give the following example.
\begin{example} \label{example-best-approx-multivariate} \rm
Let us consider the function $g_j:I \to \R$ (with $j$ fixed between 1 and $d$), defined on the d-variate rectangle $I$, as:
$$
g_j({\tt x})\ :=\ x_j-a_j, \quad \quad {\tt x} = (x_1,...,x_d) \in I.
$$
Assume for simplicity $a_i<b_i$ (defining $I$) both integers, for $i=1,...,d$. Thus we have, e.g., for ${\tt a}:=(a_1,a_2,...,a_d)\in I$:
$$
F_n(g_j,{\tt a}) - g_j({\tt a})\ =\ F_n(g_j,{\tt a})\ =\ {\displaystyle \sum_{{\tt k} \in {\cal J}_n} g_j({\tt k}/n)\,  \Psi_{\sigma}(n{\tt a}-{\tt k}) \over \displaystyle \sum_{{\tt k} \in {\cal J}_n} \Psi_{\sigma}(n{\tt a}-{\tt k})}
$$
$$
=\ {\displaystyle \sum_{k_j=na_j}^{nb_j} \left( {k_j \over n}-a_j  \right)\,  \phis(na_j-k_j) \over \displaystyle  \sum_{k_j=na_j}^{nb_j} \phis(na_j-k_j)}\ =\ n^{-1}\,	\cdot {\displaystyle \sum_{k_j=na_j}^{nb_j} \left( k_j-n a_j  \right)\,  \phis(na_j-k_j) \over \displaystyle  \sum_{k_j=na_j}^{nb_j} \phis(na_j-k_j)}.
$$
Now, setting $\nu=k_j-n a_j$ in the above sum, and recalling that $\phis$ is even and $\sum_{k_j=na_j}^{nb_j} \phis(na_j-k_j)\miu 1$, we can obtain:
$$
F_n(g_j,{\tt a}) - g_j({\tt a})\ \mau\ n^{-1}\,	\sum_{\nu=0}^{n(b_j -a_j)} \nu\,  \phis(\nu) \ \mau n^{-1}\, \phis(1)\ >\ 0
$$
from which we get:
$$
\| F_n(g_j, \cdot) - g_j(\cdot)\|_\infty \mau\ \phis(1)\,  \omega(g_j, 1/n)\ >\ 0,
$$
being $\omega(g_j, 1/n)=1/n$, $n \in \N$. The latter inequality shows that the {\em global} estimate in (\ref{order-NN}) can not be improved in general in $C(I)$, or in its sub-spaces $C^r(I)$, $r \mau 1$.

Similarly, for the Kantorovich NN operators (see \cite{COPI2}) we obtain:
$$
K_n(g_j,{\tt a}) - g_j({\tt a})\ =\ K_n(g_j,{\tt a})\ =\ {\displaystyle \sum_{{\tt k} \in {\cal K}_n} \left[ n^d\, \int_{{\tt k}/n}^{({\tt k}+1)/n} g_j({\tt u})\, d{\tt u}\right]\,  \Psi_{\sigma}(n{\tt a}-{\tt k}) \over \displaystyle \sum_{{\tt k} \in {\cal K}_n} \Psi_{\sigma}(n{\tt a}-{\tt k})}
$$
$$
=\ {\displaystyle \sum_{k_j=na_j}^{nb_j-1} \, \left[ n\, \int_{k_j/n}^{(k_j+1)/n} (u_j - a_j)\, du_j\right]\, \phis(na_j-k_j) \over \displaystyle  \sum_{k_j=na_j}^{nb_j-1} \phis(na_j-k_j)}
$$
$$
=\ n^{-1}\,	\cdot {\displaystyle \sum_{k_j=na_j}^{nb_j-1} \left[\frac12 + \left( k_j-n a_j  \right)\right]\,  \phis(na_j-k_j) \over \displaystyle  \sum_{k_j=na_j}^{nb_j-1} \phis(na_j-k_j)}\
\mau\ \omega(g_j,1/n)\, \left[ \frac12 + \phis(1)  \right]\ >\ 0,
$$
deducing the same conclusions achieved for $F_n$.
\end{example}
\begin{remark} \rm
From the examples given in Example \ref{example-best-approx-multivariate} we can deduce {\em lower bounds} for the constant $C_1$ and $C_2$ of inequalities (\ref{order-NN}) and (\ref{order-KNN}). At the same time, retracing the proofs given in \cite{COCOKA1}, we can finally obtain the following two-sided bound for the {\em constants of best approximation}: 
$$
0<\, \phis(1)\ \miu\ C_1\ \miu\ \phis(1)^{-d}\, \cdot \left[ 1\, +\, M_1(\Psi_\sigma)  \right],
$$
and 
$$
0<\, \phis(1)\ + \frac12\ \miu\ C_2\ \miu\ \phis(2)^{-d}\, \cdot \left[ 2\, +\, M_1(\Psi_\sigma)  \right].
$$
\end{remark}
%


\section{Local Voronovskaja-type theorems} \label{sec4}

We recall that the Fourier transform of an $L^1$-function $f:\R \to \R$ is defined as follows:
$$
\widehat{f}(v)\ :=\ \int_\R f(x)\, e^{-ixv}\, dx, \quad v \in \R.
$$
The following useful property of the Fourier transform of $\phis$ can now be deduced.
\begin{lemma} \label{Fourier-lemma4.1}
Let $\alpha>1$; then:
$$
\widehat{\phis}(0)\, =\, 1, \quad \quad \widehat{\phis}(2\pi k)\, =\, 0, \quad k \in \Z\setminus\{ 0\}, \quad \quad \left[\widehat{\phis}\right]'(0)\, =\, 0.
$$
\end{lemma}
\begin{proof}
By the Poisson summation formula we know that:
$$
\sum_{k \in \Z} \phis(u-k)\ =\ \sum_{k \in \Z} \widehat{\phis}(2\pi k)e^{i 2 \pi k u},
$$
where the series of the right is the Fourier series of the 1-periodic function on the left. Since from Lemma 2.2 of \cite{COSP2} we know that:
$$
\sum_{k \in \Z} \phis(u-k) = 1, \quad \quad u \in \R,
$$
it turns out that the above Fourier series must be constant and equal to one, hence it reduces only to the term $\widehat{\phis}(0)=1$ (namely, $\widehat{\phis}(2\pi k) = 0$, $k \neq 0$, for more details see \cite{BUST}). Finally, it is well-known that:
$$
\left[\widehat{\phis}\right]'(v)\, =\, (-i)\, \widehat{h_1}(v), \quad v \in \R, \quad \mbox{with} \quad h_1(u)=u\, \phis(u), \quad u \in \R, \quad h_1 \in L^1(\R).
$$
Hence, immediately follows that:
$$
\left[\widehat{\phis}\right]'(0)\, =\, (-i)\, \int_\R x\, \phis(x)\, dx\ =\ 0,
$$
since $\phis$ is even (see Lemma 2.1 (iii) of \cite{COSP2}). This completes the proof.
\end{proof}
From Lemma \ref{Fourier-lemma4.1}, and defining the following series (the so-called {\em algebraic moments} of order $\nu$):
$$
{\cal A}_\nu (\phis, u)\ :=\ \sum_{k \in \Z} \phis(u-k)\, (k-u)^{\nu}, \quad u \in \R, \quad \nu \in \N,
$$
(we already know that ${\cal A}_0 (\phis, u)=1$, $u \in \R$) we can prove what follows.
\begin{lemma} \label{Fourier-lemma4.2}
Let $\alpha>2$; if the condition:
\be \label{ipotesi-trasformata}
\left[\widehat{\phis}\right]^{(s)}(2\pi k)\ =\ 0, \quad k \in \Z \setminus \{ 0\}, \quad s=1,2,
\ee
holds, then:
$$
{\cal A}_1 (\phis, u)\ =\ 0, \quad {\cal A}_2 (\phis, u)\ =\ - \left[\widehat{\phis}\right]^{''}(0) >0, \quad u \in \R.
$$
In particular, if $\phis$ is band-limited with $supp\, \widehat{\phis} \subset [-2\pi, 2\pi]$, assumption (\ref{ipotesi-trasformata}) trivially holds.
\end{lemma}
\begin{proof}
The proof immediately follows by the application of the Poisson summation formula to the functions $h_1(u):=u\, \phis(u)$ and $h_2(u):=-u^2\, \phis(u)$, $u \in \R$.
\end{proof}
\begin{remark} \rm \label{remark-ipotesi-trasformata}
If $\sigma \in C^{m}(\R)$ satisfies $(\Sigma 3)$ for $\alpha>m \in \N$, $m>2$, and the property (\ref{ipotesi-trasformata}) holds for $s=1,...,m$ (or analogously $\phis$ is band-limited in $[-2\pi,2\pi]$), the result of Lemma \ref{Fourier-lemma4.2} can be generalized for high-order algebraic moments. More precisely, one has that all the moments of odd order $j$ turn out to be ${\cal A}_j (\phis, u)\ =\ 0$, while if $j$ is even we get:
$$
{\cal A}_j (\phis, u)\ =\ (-1)^{j/2} \left[\widehat{\phis}\right]^{(j)}(0) >0
$$
$u \in \R$, and $j=1,...,m$.
\end{remark}
\begin{remark} \rm
From Lemma \ref{Fourier-lemma4.1} and  Lemma \ref{Fourier-lemma4.2}, it turns out that the so-called assumption $(\Sigma 5)$ introduced in Section 3 of \cite{CACO1} is always satisfied for $m=2$ with the values explicitly established in Lemma \ref{Fourier-lemma4.2} when $\phis$ satisfies (\ref{ipotesi-trasformata}) or it is band-limited. For the sake of clarity, we recall that the above mentioned assumption $(\Sigma 5)$ concerns the so-called algebraic moments (in the one-dimensional case $d=1$) defined as in the equation (\ref{eq:alg truncat moments}) below with $k \in \Z$ in place of ${\tt k} \in {\cal J}_n$, which are assumed to be constant for $\nu=1,2$. As we will show below, the values of the algebraic moments are related to the problem of the {\em saturation order} of the considered operators in $C(I)$.
\end{remark}
Now, we recall the definition of the multivariate
{\em truncated algebraic moment} (t.a.m.)  of order $\nu \in \N$ of the function $\Psi_\sigma$ (see \cite{CACO2}), and given by:
\begin{equation}
m_{\nu, {\tt h}}^{n}\left(\Psi_\sigma, {\tt u}\right):=\sum_{{\tt k} \in {\cal J}_n}\Psi_\sigma\left( {\tt u} -{\tt k}\right)\left({\tt k}- {\tt u} \right)^{{\tt h}},\, \quad {\tt u} \in\mathbb{R}^{d},\label{eq:alg truncat moments}
\end{equation}
for every $n\in\mathbb{N}^{+}$, ${\tt h} \in \N^d_0$, $\nu=|{\tt h}|:=h_1 + ... + h_d$, where the above vectors-power is defined as:
$$
{\tt x}^{\tt h}\ :=\ \prod_{i=1}^{d} x_{i}^{h_{i}}, \quad \quad {\tt x} \in \R^d.
$$ 
Furthermore, we also introduce the following compact subset of $I$, namely:
\be \label{def-I-delta}
I_\delta\ :=\ \prod_{i=1}^d \left[ a_i+\delta,\, b_i-\delta \right]\ \subset I,
\ee
where $\delta>0$ is any fixed sufficiently small parameter.\vskip0.2cm

Now, we can state the following.
\begin{lemma}\label{trunc mom}
Let $\alpha>2$ and assume that condition (\ref{ipotesi-trasformata}) holds; then we have:
\be \label{num1}
m_{0, {\tt 0}}^{n}\left(\Psi_{\sigma},n {\tt x}\right)\, =\, 1\, +\, \mathfrak{T}_{0,0,n}\left({\tt x}\right),  
\ee
\be \label{num2}
m_{1, {\tt e}_i}^{n}\left(\Psi_{\sigma},n {\tt x}\right)\, =\, \mathfrak{T}_{1,i,n}\left({\tt x}\right),  
\ee
where ${\tt e}_i \in \R^d$, are the vectors with $e_i=1$ and $e_j=0$, $j\neq i$, $i=1,...,d$, 
\be \label{num3}
m_{2, {\tt h}_{i,j}}^{n}\left(\Psi_{\sigma},n {\tt x}\right)\, =\, \mathfrak{T}_{2,i,j,n}\left( {\tt x}\right),  
\ee
where ${\tt h}_{i,j} \in \R^d$, are the vectors with $h_i=h_j=1$ and $h_t=0$ with $t\neq i$ and $t \neq j$, $i,j=1,...,d$, and
\be \label{num4}
m_{2, {\tt s}_{i}}^{n}\left(\Psi_{\sigma},n {\tt x}\right)\, =\, - \left[\widehat{\phis}\right]^{''}(0)\ +\ \mathfrak{T}_{2,i,n}\left({\tt x}\right),  
\ee
where ${\tt s}_i \in \R^d$, are the vectors with $s_i=2$ and $s_j=0$, $j\neq i$, $i=1,...,d$. The above functions $\mathfrak{T}_{\nu,i,n}$ and $\mathfrak{T}_{\nu,i,j,n}$ have the property:
\begin{equation} \label{num5}
|\mathfrak{T}_{\nu,i,n}({\tt x})|=\, {\cal O}\left( n^{-\alpha-1} \right), \quad and \quad \left| \mathfrak{T}_{\nu,i,j,n} \left( {\tt x}\right) \right|\, =\, {\cal O}\left( n^{-\alpha-1} \right),
\end{equation}
uniformly for $\xx \in I_{\delta}$, as $n \to +\infty$.
\end{lemma}
\begin{proof}
Noting that:
$$
m_{0, {\tt 0}}^{n}\left(\Psi_{\sigma},n {\tt x}\right)\, =\ \prod_{i=1}^d \left[ \sum_{k_i=\lceil na_i \rceil}^{\lfloor nb_i \rfloor}  \phis(nx_i-k_i) \right]\ 
$$
$$
=\ \prod_{i=1}^d \left[ \frac12\{ \sigma(n[x_i-a_i]+na_i-\lceil na_i \rceil +1 )  + \sigma(n[x_i-a_i]+na_i-\lceil na_i \rceil) \} \right.
$$
\be \label{products-o}
\left. -\ \frac12 \{ \sigma(n[x_i-b_i]+nb_i-\lfloor nb_i \rfloor -1 )  + \sigma(n[x_i-b_i]+nb_i-\lfloor nb_i \rfloor)  \}\right]
\ee
since the above sums are telescoping. Since:
$$
n[x_i-a_i] \mau n \delta, \quad and \quad n[x_i-b_i] \miu -n\delta,
$$
for $x_i \in [a_i+\delta,b_i-\delta]$, recalling the properties of $\sigma$, assumption $(\Sigma 3)$ and computing all the products in (\ref{products-o}) we immediately get $(\ref{num1})$ and that:
$$
\left| \mathfrak{T}_{0, 0,n} \left( {\tt x}\right) \right|\, =\, {\cal O}\left( n^{-\alpha-1} \right), \quad as \quad n \to +\infty,
$$
uniformly in $I_\delta$.  Similarly, we have:
$$
m_{1, {\tt e}_i}^{n}\left(\Psi_{\sigma},n {\tt x}\right)\, =\, \left[  \sum_{k_i=\lceil na_i \rceil}^{\lfloor nb_i \rfloor} (k_i-nx_i) \phis(nx_i-k_i) \right]\, \prod_{j=1, j \neq i}^d \left[ \sum_{k_j=\lceil na_j \rceil}^{\lfloor nb_j \rfloor}  \phis(nx_j-k_j) \right]
$$
$$
=\ -\, \left[  \sum_{k_i<\lceil na_i \rceil} (k_i-nx_i) \phis(nx_i-k_i) + \sum_{k_i>\lfloor nb_i \rfloor} (k_i-nx_i) \phis(nx_i-k_i) \right]\, \prod_{j=1, j \neq i}^d \left[ \sum_{k_j=\lceil na_j \rceil}^{\lfloor nb_j \rfloor}  \phis(nx_j-k_j) \right],
$$
where the above equality comes from Lemma \ref{Fourier-lemma4.2} since ${\cal A}_1(\phis,nx_i)=0$. Now, setting:
$$
m_{1, {\tt e}_i}^{n}\left(\Psi_{\sigma},n {\tt x}\right)\, =\ \mathfrak{T}_{1,i,n}\left({\tt x}\right),
$$
using the above expression for the first order truncated moment, by assumption $(\Sigma 3)$ and since any $\sum_{k_j=\lceil na_j \rceil}^{\lfloor nb_j \rfloor}  \phis(nx_j-k_j) \miu 1$ we immediately obtain (\ref{num5}). Since the proof of (\ref{num3}) is completely analogous to the previous one, we immediately consider the case in (\ref{num4}). However, using again Lemma \ref{Fourier-lemma4.2} we can note that:
$$
m_{2, {\tt s}_i}^{n}\left(\Psi_{\sigma},n {\tt x}\right)\, 
$$
$$
=\ \left[ -\left[\widehat{\phis}\right]^{''}(0)\ -\ \sum_{k_i<\lceil na_i \rceil} (k_i-nx_i)^2 \phis(nx_i-k_i)\ -\ \sum_{k_i>\lfloor nb_i \rfloor} (k_i-nx_i)^2 \phis(nx_i-k_i) \right] \times
$$
$$
\prod_{j=1, j \neq i}^d \left[ \sum_{k_j=\lceil na_j \rceil}^{\lfloor nb_j \rfloor}  \phis(nx_j-k_j) \right],
$$
then the proof follows by the same reasoning used above, and observing that by the first part of the present lemma the term:
$$
\prod_{j=1, j \neq i}^d \left[ \sum_{k_j=\lceil na_j \rceil}^{\lfloor nb_j \rfloor}  \phis(nx_j-k_j) \right] \ \to \ 1, \quad as \quad n\to+\infty,
$$
uniformly in $I_\delta$.
\end{proof}

  Now we are able to state the following {\em local} Voronovskaja-type theorem.
\begin{theorem} \label{voronovskaja-NN}
Let $\alpha>3$ and assume that condition (\ref{ipotesi-trasformata}) holds; for $f \in C^2(I)$ and ${\tt x} \in I_\delta$ we have:
$$
\lim_{n \to +\infty} n^2\, \left[ F_n (f , {\tt x}) -  f({\tt x})  \right]\ =\ -\frac12\, \left[\widehat{\phis}\right]^{''}\!\!\!(0)\, \cdot (\Delta f)({\tt x}),
$$
where $\Delta f$ denotes the Laplacian of $f$. In particular, the above limit holds uniformly in $I_\delta$, and pointwise is valid for every ${\tt x} \in I^0$ (namely, in the interior of $I$).
\end{theorem}
\begin{proof}
Using the second order Taylor formula with Lagrange remainder at ${\tt x} \in I_\delta$ we have:
$$
f\left({\tt u}\right)=f({\tt x})\ +\ \sum_{i=1}^d {\partial f \over \partial x_i}({\tt x})\, (u_i-x_i)\, +\, \sum_{i=1}^d \sum_{\stackrel{j=1}{j \neq i}}^d {\partial^2 f \over \partial x_i \partial x_j}({\tt x})\, (u_i-x_i)\, (u_j-x_j)\, 
$$
\be  \label{taylor-lagrange}
+\ \frac12\, \sum_{i=1}^d {\partial^2 f \over \partial x_i^2}({\tt x})\, (u_i-x_i)^2\, +\ \sum_{|{\tt h}|=2}{(D^{{\tt h}}f)({\tt y}) - (D^{{\tt h}}f)({\tt x}) \over {\tt h}!}\, ({\tt u}-{\tt x})^{{\tt h}},
\ee
where ${\tt y}$ is a suitable vector belonging to the segment connecting ${\tt u}$ and ${\tt x}$, the symbol $D^{{\tt h}}f$ denotes any second order partial derivative of $f$. Using (\ref{taylor-lagrange}) with ${\tt u} = {\tt k}/n$ in the operators $F_n f$ we obtain:
$$
n^2\left[F_n(f, {\tt x})\ -\ f({\tt x})\right]\ =\ n\, {\displaystyle \sum_{i=1}^d {\partial f \over \partial x_i}({\tt x})\ \cdot m^n_{1,{\tt e}_i}(\Psi_\sigma, n {\xx}) \over  m^n_{0,{\tt 0}}(\Psi_\sigma, n {\xx})}\ +\ {\displaystyle \sum_{i=1}^d \sum_{\stackrel{j=1}{j \neq i}}^d {\partial^2 f \over \partial x_i \partial x_j}({\tt x})\,\ \cdot m^n_{2,{\tt h}_{i,j}}(\Psi_\sigma, n {\xx}) \over  m^n_{0,{\tt 0}}(\Psi_\sigma, n {\xx})}\
$$
$$
+\ {\displaystyle \frac12 \, \sum_{i=1}^d {\partial^2 f \over \partial x_i^2}({\tt x})\,\ \cdot m^n_{2,{\tt s}_{i}}(\Psi_\sigma, n {\xx}) \over  m^n_{0,{\tt 0}}(\Psi_\sigma, n {\xx})}\ +\, {\displaystyle \sum_{|{\tt h}|=2}({\tt h} !)^{-1}\sum_{{\tt k} \in {\cal J}_n}\left[(D^{{\tt h}}f)({\tt y}_{{\tt k}/n}) - (D^{{\tt h}}f)({\tt x})\right] \left( {\tt k}-n{\tt x} \right)^{{\tt h}}\Psi_\sigma(n {\tt x}-{\tt k}) \over  m^n_{0,{\tt 0}}(\Psi_\sigma, n {\xx})}
$$
$$
=:\ T_1\ +\ T_2\ +\ T_3\ +\ T_4.
$$
where the above symbols are those used in Lemma \ref{trunc mom}, and ${\tt y}_{{\tt k}/n}$ are suitable vectors belonging to the segment connecting ${\tt k}/n$ and ${\tt x}$. It is clear that by the boundedness of $f$ and its second order partial derivatives, by Lemma \ref{trunc mom} we immediately have that $T_1$ and $T_2$ tend to zero (uniformly with respect to ${\tt x} \in I_\delta$) as $n \to +\infty$. Concerning $T_3$, we have:
$$
\lim_{n \to +\infty} T_3\ =\ -\frac12\ \sum_{i=1}^d {\partial^2 f \over \partial x_i^2}({\tt x})\, \left[\widehat{\phis}\right]^{''}\!\!\!(0)\ =\ -\frac12\, \left[\widehat{\phis}\right]^{''}\!\!\!(0)\, (\Delta f)({\tt x}),
$$
uniformly with respect to ${\tt x} \in I_\delta$. Hence, in order to get the thesis it remains to show that $T_4$ converges to zero. We can write what follows using the well-known properties of the modulus of continuity:
$$
|T_4|\ \miu\ \phis(1)^{-d}\, \sum_{|{\tt h}|=2}\, \sum_{{\tt k} \in {\cal J}_n}\omega \left(D^{{\tt h}}f, \|{\tt y}_{{\tt k}/n}  - {\tt x}\|_2 \right) \left| {\tt k}-n{\tt x} \right|^{{\tt h}}\Psi_\sigma(n {\tt x}-{\tt k}) 
$$
$$
\miu\ \phis(1)^{-d}\, \sum_{|{\tt h}|=2}\, \sum_{{\tt k} \in {\cal J}_n}\omega \left(D^{{\tt h}}f, \|{\tt k}/n  - {\tt x}\|_2 \right) \left| {\tt k}-n{\tt x} \right|^{{\tt h}}\Psi_\sigma(n {\tt x}-{\tt k}) 
$$
$$
\miu\ \phis(1)^{-d}\, \sum_{|{\tt h}|=2}\, \omega(D^{{\tt h}}f, 1/n)\,  \sum_{{\tt k} \in {\cal J}_n} (1+ \|{\tt k} - n {\tt x}\|_2)\left| {\tt k}-n{\tt x} \right|^{{\tt h}}\Psi_\sigma(n {\tt x}-{\tt k})
$$
$$
\miu\ \phis(1)^{-d}\, \sum_{|{\tt h}|=2}\, \omega(D^{{\tt h}}f, 1/n)\,  \sum_{{\tt k} \in {\cal J}_n} (1+ \|{\tt k} - n {\tt x}\|_2)\,  \| {\tt k}-n{\tt x} \|_2^2\, \Psi_\sigma(n {\tt x}-{\tt k})
$$
$$
\miu\ \phis(1)^{-d}\, \widetilde{C}\ \sum_{|{\tt h}|=2}\, \omega(D^{{\tt h}}f, 1/n)\,   <\ +\infty,
$$
where $\widetilde{C}<+\infty$ since is related to a suitable combination of $M_\nu(\Psi_\sigma)$, $\nu=2,3$, that are all finite. Hence the proof follows by the properties of the modulus of continuity, noting that any $D^{{\tt h}}f$ is continuous on $I$ and that in the above sum we have a finite number of terms (depending on $d$).
\end{proof}
\begin{remark} \rm
Note that, in order to provide a local Vorovovskaja-type result for the operators $K_n f$ we have to consider the following definition for the truncated algebraic moments:
\be \label{eq:alg truncat moments kantorovich}
\widetilde{m}_{\nu, {\tt h}}^{n}\left(\Psi_\sigma, {\tt u}\right)\ :=\ \sum_{{\tt k} \in {\cal K}_n}\Psi_\sigma\left( {\tt u} -{\tt k}\right)\left({\tt k}- {\tt u} \right)^{{\tt h}},\, \quad {\tt u} \in\mathbb{R}^{d}.  
\ee
For $\widetilde{m}_{\nu, {\tt h}}^{n}\left(\Psi_\sigma, {\tt u}\right)$ we can prove the same results of Lemma \ref{trunc mom}; in particular we have: 
\be 
\widetilde{m}_{0, {\tt 0}}^{n}\left(\Psi_{\sigma},n {\tt x}\right)\, =\, 1\, +\, \mathfrak{L}_{0,0,n}\left({\tt x}\right),  
\ee
\be  
\widetilde{m}_{1, {\tt e}_i}^{n}\left(\Psi_{\sigma},n {\tt x}\right)\, =\, \mathfrak{L}_{1,i,n}\left({\tt x}\right),  
\ee
where the above functions $\mathfrak{L}_{\nu,i,n}$ are ${\cal O}(n^{-\alpha-1})$, as $n \to +\infty$, uniformly for ${\xx} \in I_{\delta}$.
\end{remark}

Thus, using the same strategy proposed in the proof of Theorem \ref{voronovskaja-NN} we obtain what follows.
\begin{theorem} \label{voronovskaja-KNN}
Let $\alpha>2$ and assume that condition (\ref{ipotesi-trasformata}) holds for $s=1$; for any $f \in C^1(I)$ and ${\tt x} \in I^0$ we have:
$$
\lim_{n \to +\infty} n\, \left[ K_n (f , {\tt x}) -  f({\tt x})  \right]\ =\ \frac12\, \sum_{i=1}^d {\partial f \over \partial x_i}({\tt x}).
$$
The above limit holds uniformly in $I_\delta$.
\end{theorem}
\begin{proof}
Using the first order Taylor formula with Lagrange remainder  in the operators $K_n f$ we have: 
$$
n\left[ (K_nf)({\tt x})-f(\xx)\right]\ =\ n\, {\displaystyle \sum_{i=1}^d {\partial f \over \partial x_i}({\tt x})\ \cdot \sum_{{\tt k} \in {\cal K}_n} \left[ n^d \int_{{\tt k}/n}^{({\tt k}+1)/n} (u_i-x_i)\, d{\tt u} \right] \Psi_{\sigma}(n{\tt x}-{\tt k})\over  \widetilde{m}^n_{0,{\tt 0}}(\Psi_\sigma, n {\xx})}\ 
$$
$$
+\ n\, {\displaystyle \sum_{i=1}^d  \sum_{{\tt k} \in {\cal K}_n} \left({\partial f \over \partial x_i}({\tt y}_{{\tt k}/n})\ \ -\ {\partial f \over \partial x_i}({\tt x})  \right) \cdot  \left[ n^d \int_{{\tt k}/n}^{({\tt k}+1)/n} (u_i-x_i)\, d{\tt u} \right] \Psi_{\sigma}(n{\tt x}-{\tt k})\over  \widetilde{m}^n_{0,{\tt 0}}(\Psi_\sigma, n {\xx})}\
$$
$$
=\  n\, {\displaystyle \sum_{i=1}^d {\partial f \over \partial x_i}({\tt x})\ \cdot \sum_{{\tt k} \in {\cal K}_n} \left[ n \int_{k_i/n}^{(k_i+1)/n} (u_i-x_i)\, du_i \right] \Psi_{\sigma}(n{\tt x}-{\tt k})\over  \widetilde{m}^n_{0,{\tt 0}}(\Psi_\sigma, n {\xx})}
$$
$$
+\ n\,  {\displaystyle \sum_{i=1}^d  \sum_{{\tt k} \in {\cal K}_n} \left({\partial f \over \partial x_i}({\tt y}_{{\tt k}/n})\ \ -\ {\partial f \over \partial x_i}({\tt x})  \right) \cdot  \left[ n \int_{k_i/n}^{(k_i+1)/n} (u_i-x_i)\, du_i \right] \Psi_{\sigma}(n{\tt x}-{\tt k})\over  \widetilde{m}^n_{0,{\tt 0}}(\Psi_\sigma, n {\xx})}\
$$
$$
=:\ {\displaystyle \sum_{i=1}^d {\partial f \over \partial x_i}({\tt x})\ \cdot \sum_{{\tt k} \in {\cal K}_n} \left[ \frac12 + (k_i-nx_i) \right] \Psi_{\sigma}(n{\tt x}-{\tt k})\over  \widetilde{m}^n_{0,{\tt 0}}(\Psi_\sigma, n {\xx})}\ +\ R^1_n({\xx})
$$
$$
=\ \frac12\, \sum_{i=1}^d {\partial f \over \partial x_i}({\tt x})\ +\ {\displaystyle \sum_{i=1}^d {\partial f \over \partial x_i}({\tt x})\ \cdot \widetilde{m}^n_{1, {\tt e}_i}(\Psi_\sigma, n{\tt x})   \over  \widetilde{m}^n_{0,{\tt 0}}(\Psi_\sigma, n {\xx})}\ +\ R^1_n({\xx}),
$$
where the vectors ${\tt e}_i$ are those defined in the statement of Lemma \ref{trunc mom}. Now, the proof can be easily completed following the same reasoning used in the proof of Theorem \ref{voronovskaja-NN}.
\end{proof}

Note that, the order of approximation deduced by the first-order Voronovskaja theorem for the operators $K_n$ does not improve even when we consider $f \in C^2(I)$, due to the first order term which, in general, does not converge to zero.


\section{Saturation theorems by the solution of the Laplace equation} \label{sec5}

From the results of Section \ref{sec3} and Section \ref{sec4} concerning the operators $F_n$ we know that for $n \to +\infty$, the {\em global} best possible  order of approximation in $C(I)$ is ${\cal O}(n^{-1})$, that in general can not be improved. However, restricting our analysis {\em locally} on $I_\delta$, we know that also the order ${\cal O}(n^{-2})$ can be achieved for a certain class of functions.

The above reasoning leads to the problem of establishing the {\em saturation order} of approximation for the above operators and the corresponding class of functions in which the operators are saturated. 

The latter problem can be in general difficult. In the literature, in order to approach this kind of problems one can find the so-called {\em parabola technique} introduced by Bajsanski and Bojanic in \cite{BABO1964} in order to study the saturation order for the Bernstein polynomials; later this method has been generalized to face a wide family of (univariate) linear operators (see Chapter 4 of \cite{Devore1972}). The extension of the above technique to the multivariate case has been given by Ditzian (see \cite{Ditzian1983}) and it is known as the {\em generalized parabola technique}. On the basis of the latter approach there is a useful theorem (see Theorem 2.1 of \cite{Ditzian1983}); the proof of such a result is based on the fact that the domain of $f$ is a ball with respect to the classical Euclidean norm of $\R^d$. Since a multivariate rectangle as $I$, topologically speaking, is substantially equivalent to a ball, it is not difficult to show that the result of Ditzian is still valid for functions $f:I\to \R$ as those considered in this paper. However, from the proof in this context we can see that the generalized parabola function can be written in a certain specific way that can turn out to be useful in certain situations. 

For the latter reason it can be useful to explicitly state the proof of what follows.
\begin{theorem} \label{th-Ditzian-on-I}
Let $f \in C(I)$, such that $f({\tt x})=0$ for ${\tt x} \in \partial I$ (the boundary of $I$) and $f({\tt x}_0)>0$ for a certain ${\tt x}_0$ belonging to $I^{0}$. Then, there exists a function 
\be \label{generalized-parabola2}
Q({\tt x})\ :=\ A\, \sum_{i=1}^d(x_i-a_i)^2\ +\ C \ =:\ A\, \left(  \sum_{i=1}^d x_i^2  \right)\, +\, \sum_{i=1}^d B_i\, x_i\ +\ \bar \gamma,
\ee
with $A<0$, $C \in \R$, $Q({\tt x})>0$ on $I$, $Q({\tt x}) \mau f({\tt x})$ on $I$ and there exists ${\tt y} \in I^0$ such that $Q({\tt y})=f({\tt y})$.
\end{theorem}
\begin{proof}
Let $M>0$ be the maximum of $f$ on $I$. We now define:
$$
Q_1({\tt x})\ :=\ A\, \sum_{i=1}^d(x_i-a_i)^2\ +\ M\ - A\, r^2,
$$
where $r>0$ is chosen as the radius of a ball centered at ${\tt a}=(a_1,...,a_d)$ and such that contains $I$. Moreover, we also fix $A<0$ sufficiently small such that $-Ar^2 < f({\tt x}_0)$. By construction we have $Q_1({\tt x}) \mau M$, since we chosen $r$ such that $\max_{{\tt x} \in I}\sum_{i=1}^d(x_i-a_i)^2 \miu r^2$. Let now:
$$
m\ :=\ \min_{{\tt x} \in I}\left[ Q_1({\tt x}) - f({\tt x}) \right],
$$
and we define:
$$
Q({\tt x})\ :=\ Q_1({\tt x})\, -\, m\ =:\ A\, \sum_{i=1}^d(x_i-a_i)^2\ +\ C, \quad {\tt x} \in I,
$$
with $C:=M- A r^2-m$. Obviously, we have: $Q({\tt x}) \mau f({\tt x})$, ${\tt x} \in I$, by the definition of $m$. Moreover, for any ${\tt x} \in \partial I$, recalling that $f({\tt x})=0$ and $Q_1({\tt x}) \mau M$, we also have $Q_1({\tt x})-f({\tt x}) \mau M$. However, we also know that in ${\tt x}_0=(x^0_1,...,x^0_d) \in I^0$:
\be \label{per-il-minimo}
Q_1({\tt x}_0)-f({\tt x}_0)\ =\ A\, \sum_{i=1}^d(x^0_i-a_i)^2\ +\ M\ - A\, r^2 -\, f({\tt x}_0)\ <\ M\ - A\, r^2 -\, f({\tt x}_0)\ <\ M.
\ee
The above considerations shown that the minimum $m$ of the function $Q_1({\tt x})-f({\tt x})$ must be achieved in a point ${\tt y} \in I^0$, namely:
$$
Q_1({\tt y})-f({\tt y})\ =\ m \quad \mbox{that implies} \quad Q({\tt y}) = Q_1({\tt y})- m\ = f({\tt y}).
$$
At the same time (\ref{per-il-minimo}) implies that the minimum $m<M$, and since $Q_1({\tt x}) \mau M$ we also deduce $Q({\tt x})>0$ on $I$. This completes the proof.
\end{proof}
   We are now able to prove the following.
\begin{theorem} \label{th-saturation-NN-local}
Let $\alpha>3$, assume that condition (\ref{ipotesi-trasformata}) holds, and let $f \in C(I)$ be fixed. If:
\be \label{ipotesi-o-small-locale}
\sup_{{\tt x} \in I_\delta} |F_n(f, {\tt x} )-f({\tt x})|\ =\ o(n^{-2}), \quad as \quad n \to +\infty,
\ee
then $f \in C^2(I_\delta)$ and is harmonic in $I_\delta$, i.e., $f$ is a classical solution of the Laplace equation on $I_\delta$:
$$
(\Delta f)({\tt x})\, =\, 0, \quad {\tt x} \in I_\delta.
$$
Vice-versa, if $f \in C^2(I_\delta)$ is harmonic, then (\ref{ipotesi-o-small-locale}) holds.
\end{theorem}
\begin{proof}
First, we can note that if the function $f$ is the classical solution of the Laplace equation, by Theorem \ref{voronovskaja-NN} we immediately get the second part of the thesis.
\vskip0.2cm

We now prove the first part of the thesis.
Consider the $C^2$-function $g: I_\delta \to \R$ which is the classical solution of the following Dirichlet problem for the Laplace equation on $I_\delta$: 
$$
\begin{cases}
(\Delta g)({\tt x})\, =\, 0, \quad {\tt x} \in I_\delta^0, \\
g({\tt x})\, =\, f({\tt x}), \quad {\tt x} \in \partial I_\delta,
\end{cases}
$$
with continuous boundary condition $f$ on $\partial I_\delta$. We now set $G({\tt x}):=f({\tt x})-g({\tt x})$, ${\tt x} \in I_\delta$; clearly it turns out that $G \in C(I_\delta)$. By the above construction, we have that $G({\tt x})=0$ on $\partial I_\delta$, and this implies that $G$ assumes in $I_\delta^0$ its maximum or its minimum, otherwise it would be $G\equiv 0$ on $I_\delta$ and the thesis would be trivially satisfied. 

  Suppose that $G$ assumes its maximum in ${\tt x}_0 \in I_\delta^0$ with $G({\tt x}_0)>0$ (the proof will be analogous if ${\tt x}_0$ is a minimum point, replacing in the computations below $G$ with $-G$). 

  Applying Theorem \ref{th-Ditzian-on-I} to $G$ on $I_{\delta}$, we know that there exists a function
$$
Q({\tt x})\ :=\ A\, \left(  \sum_{i=1}^d x_i^2  \right)\, +\, \sum_{i=1}^d B_i\, x_i\ +\ {\bar \gamma},
$$
with $A<0$, $Q({\tt x})>0$, $Q({\tt x}) \mau G({\tt x})$ on $I_{\delta}$ and such that there exists ${\tt y} \in I_{\delta}^0$ with $Q({\tt y})=G({\tt y})$.

   Let $\gamma>0$ sufficiently small such that ${\tt y} \in I_{\delta + \gamma}$, where the set $I_{\delta + \gamma}$ is defined according to (\ref{def-I-delta}). Now, to avoid confusion, we denote by $\bar F_n$ the NN operators $F_n$ considered on the basic set $I_\delta$. For instance, for $g$ we can write:
$$
\bar F_n(g,{\tt x})\ =\ {\displaystyle \sum_{{\tt k} \in \widetilde {\cal J}_n} g({\tt k}/n)\, \Psi_\sigma(n{\tt x}-{\tt k}) \over \displaystyle \sum_{{\tt k} \in \widetilde {\cal J}_n} \Psi_\sigma(n{\tt x}-{\tt k})}, \quad {\tt x} \in I_\delta,
$$
where here:
$$
\widetilde {\cal J}_n\ :=\ \left\{ {\tt k} \in \Z^d:\ \lceil n(a_i+\delta)\rceil \miu k_i \miu \lfloor n(b_i-\delta)\rfloor ,\ i=1,...,d \right\}.
$$
As a consequence of Theorem \ref{voronovskaja-NN}, since $g \in C^2(I_\delta)$ and $(\Delta g)({\tt x})= 0$ on $I_{\delta+\gamma} \subset I_\delta$, we have:
\be \label{o-piccolo-F-bar}
\sup_{{\tt x} \in I_{\delta+\gamma}} \left| \bar F_n(g,{\tt x}) - g({\tt x})\right|\ =\ o(n^{-2}), \quad as \quad n \to +\infty,
\ee
Now, we claim that also:
\be \label{claim-ausiliare}
\sup_{{\tt x} \in I_{\delta+\gamma}} \left| \bar F_n(f,{\tt x}) - f({\tt x})\right|\ =\ o(n^{-2}), \quad as \quad n \to +\infty.
\ee
Indeed, for ${\tt x} \in I_{\delta+\gamma}$ we can write what follows:
\be \label{ineq-bar-F-1}
\bar F_n(f,{\tt x}) - f({\tt x})\ =\ {\displaystyle \sum_{{\tt k} \in {\cal J}_n} \left[f({\tt k}/n) - f({\tt x}) \right]\Psi_\sigma(n{\tt x}-{\tt k})\ -\  \sum_{{\tt k} \in {\cal J}_n \setminus \widetilde {\cal J}_n} \left[f({\tt k}/n) - f({\tt x}) \right]\Psi_\sigma(n{\tt x}-{\tt k}) \over \displaystyle \sum_{{\tt k} \in \widetilde {\cal J}_n} \Psi_\sigma(n{\tt x}-{\tt k})},
\ee
where ${\cal J}_n$ is the usual set of indexes defined considering $I$ as the domain of $f$. Thus:
$$
\hskip-0.6cm n^{2}\, \left|\bar F_n(f,{\tt x}) - f({\tt x})\right|\ \miu\ n^{2}\,  \left|   {\displaystyle   \sum_{{\tt k} \in {\cal J}_n} \left[f({\tt k}/n) - f({\tt x}) \right]\Psi_\sigma(n{\tt x}-{\tt k})  \over \displaystyle \sum_{{\tt k} \in {\cal J}_n} \Psi_\sigma(n{\tt x}-{\tt k})} \right|
 \cdot {\displaystyle \sum_{{\tt k} \in {\cal J}_n} \Psi_\sigma(n{\tt x}-{\tt k})  \over \displaystyle \sum_{{\tt k} \in \widetilde {\cal J}_n} \Psi_\sigma(n{\tt x}-{\tt k}) }    
$$
\be \label{ineq-bar-F-2}
+\ n^{2}\,  \left|  {\displaystyle \sum_{{\tt k} \in {\cal J}_n \setminus \widetilde {\cal J}_n} \left[f({\tt k}/n) - f({\tt x}) \right]\Psi_\sigma(n{\tt x}-{\tt k}) \over \displaystyle \sum_{{\tt k} \in \widetilde {\cal J}_n} \Psi_\sigma(n{\tt x}-{\tt k})}    \right|
\ee
$$
\miu\ {n^{2} \over \phis(1)^d}\, \left|F_n(f,{\tt x}) - f({\tt x})\right|\ +\ n^{2}\,  \left|  {\displaystyle \sum_{{\tt k} \in {\cal J}_n \setminus \widetilde {\cal J}_n} \left[f({\tt k}/n) - f({\tt x}) \right]\Psi_\sigma(n{\tt x}-{\tt k}) \over \displaystyle \sum_{{\tt k} \in \widetilde {\cal J}_n} \Psi_\sigma(n{\tt x}-{\tt k})}    \right|\ =:\ S_1\ +\ S_2.
$$
Now, if we denote respectively by $N_1$ and $N_2$ the minimum and the maximum of the set $\widetilde {\cal J}_n$, for $S_2$ we have:
$$
S_2\, \miu\ 2\, \|f\|_\infty\, \phis(1)^{-d}\, n^{2}\,  \sum_{{\tt k} \in {\cal J}_n \setminus \widetilde {\cal J}_n} \Psi_\sigma(n{\tt x}-{\tt k})
$$
$$
\miu\ 2\, \|f\|_\infty\, \phis(1)^{-d}\, n^{2}\,  \left\{  \sum_{i=1}^d \left( \sum_{k_i=\lceil na_i\rceil}^{N_1-1} \phis(nx_i-k_i) \sum_{{\tt k}_{[i]} \in {\cal J}_n^{[i]}} \Psi_\sigma^{[i]}(n{\tt x}_{[i]}-{\tt k_{[i]}}) \right)\ \right.
$$
$$
\left. +\   \sum_{i=1}^d \left( \sum_{k_i= N_2+1}^{\lfloor nb_i\rfloor } \phis(nx_i-k_i) \sum_{{\tt k}_{[i]} \in {\cal J}_n^{[i]}} \Psi_\sigma^{[i]}(n{\tt x}_{[i]}-{\tt k_{[i]}}) \right)\ \right\}
$$
$$
\miu\ 2\, \|f\|_\infty\, \phis(1)^{-d}\, n^{2}\,  \left\{  \sum_{i=1}^d \left( \sum_{k_i=\lceil na_i\rceil}^{ N_1-1 } \phis(nx_i-k_i)  \right)\ +\   \sum_{i=1}^d \left( \sum_{k_i= N_2+1}^{\lfloor nb_i\rfloor } \phis(nx_i-k_i) \right)\ \right\}.
$$
Observe now that $N_1-1$ coincides with $\lfloor n(a_i+\delta)\rfloor$ or $\lfloor n(a_i+\delta)\rfloor - 1$; thus in the first case, 
recalling the definition of $\phis$, we can observe what follows:
$$
\sum_{k_i=\lceil na_i\rceil}^{\lfloor n(a_i+\delta)\rfloor } \phis(nx_i-k_i)\ =\ \frac12\, \left[ \sigma(nx_i-\lceil na_i\rceil) - \sigma(nx_i-\lfloor n(a_i+\delta)\rfloor )    \right]
$$
$$
+\ \frac12\left[\sigma(nx_i-\lceil na_i\rceil+1)  - \sigma(nx_i-\lfloor n(a_i+\delta)\rfloor -1)  \right]\ =\ {\cal O}(n^{-\alpha-1}), \quad as \quad n \to +\infty,
$$
uniformly with respect to ${\tt x} \in I_{\delta+\gamma}$, for every $i=1,...,d$, since 
$$
nx_i-\lceil na_i\rceil + 1\ \mau\ nx_i-\lceil na_i\rceil\ =\ n(x_i - a_i)+na_i-\lceil na_i\rceil\ \mau\ n(\delta + \gamma ) - 1,
$$
$$
nx_i-\lfloor n(a_i+\delta)\rfloor \mau nx_i-\lfloor n(a_i+\delta)\rfloor - 1 \mau n \gamma + n(a_i+\delta) - \lfloor n(a_i+\delta)\rfloor - 1\ \mau n \gamma-1,
$$
and (\ref{ordine-simmetrico}) holds. Note that, if $N_1-1$ coincides with $\lfloor n(a_i+\delta)\rfloor - 1$ we can obtain analogous considerations. Similarly, $N_2+1$ can be $\lceil n(b_i-\delta)\rceil$ or $\lceil n(b_i-\delta)\rceil +1$; thus, in the first case we obtain:
$$
\sum_{k_i=\lceil n(b_i-\delta)\rceil}^{\lfloor nb_i\rfloor } \phis(nx_i-k_i)\ =\ \frac12 \left[ \sigma(nx_i- \lceil n(b_i-\delta)\rceil) + \sigma(nx_i- \lceil n(b_i-\delta)\rceil +1)\right]
$$
$$
-\ \frac12\left[ \sigma(nx_i-  \lfloor nb_i\rfloor  ) + \sigma(nx_i- \lfloor nb_i\rfloor    -1)  \right]\ =\ {\cal O}(n^{-\alpha-1}), \quad as \quad n \to +\infty,
$$
uniformly with respect to ${\tt x} \in I_{\delta+\gamma}$, for every $i=1,...,d$, since  
$$
nx_i-  \lfloor nb_i\rfloor - 1\ \miu\ nx_i-  \lfloor nb_i\rfloor \ =\ n(x_i-b_i) + nb_i -  \lfloor nb_i\rfloor\ \miu\ -n (\delta +\gamma)\, +\, 1,
$$
$$
nx_i- \lceil n(b_i-\delta)\rceil\ \miu\ nx_i- \lceil n(b_i-\delta)\rceil + 1\ \miu\ n(x_i - b_i + \delta ) + n(b_i-\delta) - \lceil n(b_i-\delta)\rceil + 1\  \miu\ -n \gamma+1, 
$$
and $(\Sigma 3)$ holds for $\alpha>3$. Again, analogous considerations can be done if $N_2+1$ is $\lceil n(b_i-\delta)\rceil +1$. The above inequalities show that $S_2={\cal O}(n^{-\alpha+1})$, uniformly on $I_{\delta+\gamma}$, as $n \to +\infty$. Noting that by (\ref{ipotesi-o-small-locale}) also $S_1 \to 0$, uniformly on $I_{\delta+\gamma}$, as $n \to +\infty$, we can finally deduce (\ref{claim-ausiliare}). As a consequence of (\ref{o-piccolo-F-bar}) and (\ref{claim-ausiliare}) we also obtain that:
\be \label{per-la-tesi-finale}
\sup_{{\tt x} \in I_{\delta+\gamma}} \left|  \bar F_n(G,  {\xx})   - G({\xx}) \right|\ =\ o(n^{-2}), \quad as \quad n \to +\infty.
\ee
Now we can write what follows:
$$
\bar F_n(G,{\tt y}) - G({\tt y})\ =\ \bar F_n(G,{\tt y}) - Q({\tt y})\ \miu\ \bar F_n(Q,{\tt y}) - Q({\tt y}), 
$$
since $G({\tt x})\miu Q({\tt x})$, on $I_\delta$ and $\bar F_n$ is a positive linear operator. Setting:
$$
E_{j,i}({\tt x})\ :=\ (x_i)^j, \quad {\tt x} \in \R^d, \quad i=1,..., d, \quad j=1,2,
$$
and using both the linearity of $\bar F_n$ and the definition of $Q$ we obtain:
$$
n^2\, \left[  \bar F_n(G,{\tt y}) - G({\tt y})\  \right] \miu\ n^2\, \left[ \bar F_n(Q,{\tt y}) - Q({\tt y}) \right]
$$
\be \label{inequality-for-contradiction}
=\ A\, \sum_{i=1}^d n^2\, \left[ \bar F_n(E_{2,i}, {\tt y}) - E_{2,i}({\tt y}) \right]\ +\ \sum_{i=1}^d B_i\, n^2\, \left[ \bar F_n(E_{1,i}, {\tt y}) - E_{1,i}({\tt y}) \right].
\ee
Observing that, $E_{j,i} \in C^2(I_\delta)$ with $(\Delta E_{1,i})({\tt y})=0$, $(\Delta E_{2,i})({\tt y})={\partial^2 E_{2,i} \over \partial x_i^2}({\tt y})=2$, ${\tt y} \in I_{\delta+\gamma}$, passing to the limit for $n \to +\infty$ in (\ref{inequality-for-contradiction}), using Theorem \ref{voronovskaja-NN} again and (\ref{per-la-tesi-finale}), we obtain:
$$
0\ \miu\ -A\, d\, \left[  \widehat{\phis} \right]''\!\!(0)\ < 0,
$$
since $A<0$ and $-\left[  \widehat{\phis} \right]''\!\!(0)>0$ by Lemma \ref{Fourier-lemma4.2}, hence we get a contradiction which implies that $G=0$ in $I_\delta$, i.e., $f=g$. This completes the proof.
\end{proof}

From Theorem \ref{th-saturation-NN-local} we deduce that the operators $F_n$ are saturated when harmonic functions are considered. 
\begin{remark} \rm
Based on the latter result, it can be natural to ask if such a saturation theorem can be extended to the whole domain of $f$. The  answer to this question is, of course, affirmative if we are in the case $d=1$. Indeed, if $f \in C(I)$ and satisfies (\ref{ipotesi-o-small-locale}) globally on the domain $I$ (in place of a subinterval), for every fixed $\delta>0$, by Theorem \ref{th-saturation-NN-local} we obtain that $f \in C^2$ with $f''(x)=0$, namely $f(x)=mx+q$, on $I_\delta$. Since $\delta$ is arbitrary and by the continuity of $f$, it turns out that $f(x)=mx+q$ on $I$. However, if $m \neq 0$ we can repeat the same computations done in Example \ref{example-best-approx-multivariate} obtaining:
$$
(F_nf)(a)-f(a)\, =\, {m \over n} {\displaystyle \sum_{k \in {\cal J}_n} (k-na)\, \phis(na-k) \over \displaystyle \sum_{k \in {\cal J}_n} \phis(na-k)}\, \mau\ {m \over n}\, {\displaystyle \sum_{k \in {\cal J}_n} (k-\lceil na \rceil )\, \phis(na-k) \over \displaystyle \sum_{k \in {\cal J}_n} \phis(na-k)}
$$
$$
\mau\ {m \over n}\, \sum_{k \in {\cal J}_n} (k-\lceil na \rceil )\, \phis(\lceil na \rceil -k) \mau\, m\, n^{-1}\, \phis(1)\ >\ 0,
$$
if $m>0$, or similarly
$$
f(a) - (F_nf)(a)\, =\, {-m \over n} {\displaystyle \sum_{k \in {\cal J}_n} (k-na)\, \phis(na-k) \over \displaystyle \sum_{k \in {\cal J}_n} \phis(na-k)}\, \mau\, -m\, n^{-1}\, \phis(1)>0,
$$
if $m<0$, for every $n \in \N$, which contradict the global version of assumption (\ref{ipotesi-o-small-locale}). This proves that $f$ must be constant.
\end{remark}
Based on what has been observed in the above remark, we can state the following corollary.
\begin{corollary} \label{cor-saturation-nuovo-d-1}
Let $\alpha>3$ and $d=1$. Moreover, assume that condition (\ref{ipotesi-trasformata}) holds, and let $f \in C(I)$ be fixed. If:
\be \label{ip-nuova-saturazione-globale}
\sup_{x \in I} |F_n(f, x)-f(x)|\ =\ o(n^{-2}), \quad as \quad n \to +\infty,
\ee
then $f$ is constant on $I=[a,b]$. Obviously, for constant functions $f$, (\ref{ip-nuova-saturazione-globale}) is trivially satisfied.
\end{corollary}

  We now recall that (according to \cite{COHI1962,LOSC1972}) a function $f \in C^2(I^0) \cap C(I)$ is called {\em sub-harmonic} if for each $I_\delta \subset I$ one has $f({\tt x}) \miu \phi({\tt x})$, ${\tt x} \in I_\delta^0$, where $\phi$ is the solution of the following Dirichlet problem:
\be \label{Dirichlet-problem-phi}
\begin{cases}
(\Delta \phi)({\tt x})\ =\ 0,\ \quad {\tt x} \in I_\delta^0, \\
\phi({\tt x})\ =\ f({\tt x}),\ \quad \hskip0.1cm {\tt x} \in \partial I_\delta.
\end{cases}
\ee
The sub-harmonic functions are characterized by the property $(\Delta f)({\tt x}) \mau 0$, for ${\tt x} \in I^0$. We can now prove the following.
\begin{theorem} \label{local-inverse-NN-sunharmonic}
Let $\alpha>3$, assume that condition (\ref{ipotesi-trasformata}) holds, and let $f \in C^2(I^0) \cap C(I)$ be fixed. If there exists $n_0 \in \N$ such that, for every ${\tt x} \in I^0$:
\be
n^2\, \left[ F_n(f,{\tt x}) - f({\tt x})  \right]\ \mau\ -M, \quad \mbox{for every } n \mau n_0,
\ee
with a fixed $M>0$, then the function:
$$
G({\tt x})\ :=\ {-2\, M \over \phis(1)^d\, [\widehat{\phis}]''(0)}u({\tt x})\, +\, f({\tt x})\ =:\ C_{\phis}\, u({\tt x})\, +\, f({\tt x}), 
$$
is sub-harmonic, where $u$ is any fixed classical solution of the elliptic equation $(\Delta u)({\tt x})=1$.
\end{theorem}
\begin{proof}
Suppose by contradiction that the function $G$ is not sub-harmonic. Then, there exists a set $I_\delta$ for which we have $G({\tt x}_0)>\phi({\tt x}_0)$, for a certain ${\tt x}_0 \in I_\delta^0$, where $\phi$ solves (\ref{Dirichlet-problem-phi}) with boundary value in $\partial I_\delta$ equal to $G$. Now, consider the auxiliary function $g({\tt x}):=G({\tt x})-\phi({\tt x})$, defined on $I_\delta$. It is clear that $g$ satisfies the assumptions of Theorem \ref{th-Ditzian-on-I}, hence there exists a function $Q({\tt x})$ of the form (\ref{generalized-parabola2}), such that $Q > 0$ on $I_\delta$, $Q({\tt x}) \mau g({\tt x})$, and $Q({\tt y}) = g({\tt y})$, for a suitable ${\tt y} \in I_\delta^0$. 

Let now $\gamma>0$ such that ${\tt y} \in I_{\delta+\gamma} \subset I_\delta$, and denote by $\bar F_n f$ the NN operators on $I_\delta$ (as made in the proof of Theorem \ref{th-saturation-NN-local}). Proceeding as in (\ref{ineq-bar-F-1}) and (\ref{ineq-bar-F-2}), we can immediately see that:
$$
n^2\left[  \bar F_n (f, {\tt x})-f({\tt x})\right]\, \mau\, -\, M\, {\displaystyle \sum_{{\tt k} \in {\cal J}_n} \Psi_\sigma(n{\tt x}-{\tt k})  \over \displaystyle \sum_{{\tt k} \in \widetilde {\cal J}_n} \Psi_\sigma(n{\tt x}-{\tt k}) }\  +\ o(1)\ \mau\ - {M \over [\phis(1)]^d}\, +\ o(1),
$$
for every ${\tt x} \in I_{\delta+\gamma}$. Hence, we can write what follows:
$$
- {M \over [\phis(1)]^d}\, +\ o(1)\ \miu\  n^2\left[  \bar F_n (f, {\tt y})-f({\tt y})\right]\, =\ n^2\left[  \bar F_n (g, {\tt y})-g({\tt y})\right]\ 
$$
$$
+\ n^2\left\{  \bar F_n (\phi-C_{\phis} u, {\tt y}) -\left[\phi({\tt y})-C_{\phis} u({\tt y}) \right]\right\}
$$
$$
\miu\ n^2\left[  \bar F_n (Q, {\tt y})-Q({\tt y})\right]\ +\ n^2\left\{  \bar F_n (\phi-C_{\phis} u, {\tt y}) -\left[\phi({\tt y})-C_{\phis} u({\tt y}) \right]\right\}.
$$
Passing to the limit as $n\to+\infty$ in the above inequality, and using Theorem \ref{voronovskaja-NN} we get:
$$
- {M \over [\phis(1)]^d}\, \miu\ -\left[ \widehat{\phis}\right]''(0)\, \left[A\, d\, +\, \frac12\, \Delta\left(\phi - C_{\phis}\, u\right)   \right]\ 
$$
$$
=\ -\left[ \widehat{\phis}\right]''(0)\, \left[A\, d\, -\, \frac12\, C_{\phis}\, \right]\ =\ - A\, d\, \left[ \widehat{\phis}\right]''(0)\, -\, {M \over [\phis(1)]^d},
$$
which gives a contradiction. This completes the proof.
\end{proof}
From Theorem \ref{local-inverse-NN-sunharmonic} we can easily deduce the following characterization.
\begin{theorem} \label{th-5.6-inv}
Let $\alpha>3$, assume that condition (\ref{ipotesi-trasformata}) holds,  and let $f \in C^2(I^0) \cap C(I)$ be fixed. If there exists $n_0 \in \N$ such that, for every ${\tt x} \in I^0$:
\be \label{disg-local-lip2}
n^2\, \left| F_n(f,{\tt x}) - f({\tt x})  \right|\ \miu\ M, \quad \mbox{for every } n \mau n_0,
\ee
then the functions
$$
G({\tt x})\ =\ C_{\phis}\, u({\tt x})\, +\, f({\tt x}), \quad and \quad L({\tt x})\ =\ C_{\phis}\, u({\tt x})\, -\, f({\tt x}), \quad C_{\phis}:={-2\, M \over \phis(1)^d\, [\widehat{\phis}]''(0)},
$$
are sub-harmonic, where $u$ is any fixed classical solution of the elliptic equation $(\Delta u)({\tt x})=1$.

   Conversely, if the functions $G$ and $L$ (defined by any general positive constant denoted again by $C_{\phis}$) are sub-harmonic, then (\ref{disg-local-lip2}) holds with $M=-2 \left[ \widehat{\phis}\right]''(0)\, C_{\phis}$.
\end{theorem}
\begin{proof}
The necessary condition immediately follows by Theorem \ref{local-inverse-NN-sunharmonic}. Concerning the sufficient condition, since $G$ and $L$ are sub-harmonic we have:
$$
\Delta\left( C_{\phis}\, u\, \pm\, f  \right)\ =\ C_{\phis}\, \pm\ \Delta\left(f  \right) \mau 0,
$$
namely $|\Delta f |\miu C_{\phis}$. Hence, using (pointwise) Theorem \ref{voronovskaja-NN} we immediately get (\ref{disg-local-lip2}).
\end{proof}
In the special case of $d=1$, in the above theorem the $C^2$-requirement on $f$ can be removed. To prove this result, we need to recall the notion of the second order Lipschitz classes defined by means of $\omega_2(f,\gamma)$, i.e., the second order modulus of smoothness of $f$, that here is directly recalled in the one-dimensional case. More precisely, we have:
$$
\omega_2(f, \gamma)\, :=\, \sup\{ |f(x+h)- 2f(x)+f(x-h)|:\ x-h,\ x,\ x+h \in I,\ and\ 0<h \miu \gamma \}, \quad \gamma>0,
$$ 
and the corresponding Lipschitz classes (see \cite{DELO1}):
$$
Lip(\nu)\ :=\ \left\{ f \in C(I):\, \omega_2(f, \gamma)\ =\ {\cal O}(\gamma^{\nu}), \ as \ \gamma \to 0^+ \right\}, \quad 0<\nu\miu 2,
$$
namely, there exists $L>0$ (the so-called {\em Lipschitz constant) and $\bar \gamma>0$ such that:
\be
\omega_2(f, \gamma)\ \miu L\, \gamma^\nu, \quad 0<\gamma  < \bar \gamma.
\ee
}
The following inverse theorem in the case corresponding to the (local) saturation order can be stated. 
\begin{theorem} \label{saturation-order-d=1}
Let $d=1$, $\alpha>3$ and assume that condition (\ref{ipotesi-trasformata}) holds.  For any $f \in C(I)$, with $I=[a,b]$, if there exists $n_0 \in \N$ such that:
\be \label{ordine-locale-su-aperto}
n^2\, \left| F_n(f,x)-f(x) \right|\ \miu\ M, \quad \mbox{for every } \quad a<x<b, \quad \mbox{and} \quad n\mau n_0,
\ee
then functions $t_{\pm}(x):=\frac12 C_{\phis}e_2(x) \pm f(x)$, are convex, where $e_2(x):=x^2$ and $C_{\phis}$ is the constant defined in the statement of Theorem \ref{th-5.6-inv}. Furthermore, it turns out that $f \in Lip(2)$ on $(a,b)$, with Lipschitz constant equal to $C_{\phis}$.
\end{theorem}
\begin{proof}
Suppose by contradiction that $t_+(x)=\frac12 C_{\phis}e_2(x) + f(x)$ is not convex. 
Then, there exists points $a<x_1<y<x_2<b$ such that the linear function $\ell(x)$ connecting $x_1$ and $x_2$ satisfies $\ell(y)<t_+(y)$. Let now $g(x):=t_+(x)-\ell(x)$ the continuous function defined in $[x_1,x_2]$, vanishing at the endpoints and with $g(y)>0$. Using again the (generalized) parabola method in the one-dimensional case, we know that there exists a non-negative $Q(x)=A(x-x_1)^2+C$, with $A<0$, such that $g(x)\miu Q(x)$, and $g(z)=Q(z)$, for a certain $x_1 < z < x_2$. Proceeding as in the proof of Theorem \ref{local-inverse-NN-sunharmonic}, and considering the operators $\bar F_n$ restricted on $[x_1,x_2]$ and the Voronovskaja-type theorem we get:
$$
o(1)\, -\, {M \over \phis(1)}\ \miu\ n^2\, \left[ \bar F_n(f, z) - f(z)   \right]\ =\  n^2\, \left[ \bar F_n(t_+, z) - t_+(z)   \right] - \frac{n^2}{2}\, \left[ \bar F_n(C_{\phis}\, e_2, z) - C_{\phis}\, e_2(z) \right]
$$
$$
=\ n^2\, \left[ \bar F_n(g, z) - g(z)   \right]\, +\, n^2\, \left[ \bar F_n(\ell, z) - \ell(z)   \right]\, - \frac{n^2}{2}\, \left[ \bar F_n(C_{\phis}\, e_2, z) - C_{\phis}\, e_2(z) \right]
$$
$$
\miu\ n^2\, \left[ \bar F_n(Q, z) - Q(z)   \right]\, +\, n^2\, \left[ \bar F_n(\ell, z) - \ell(z)   \right]\, - n^2\, \left[ \bar F_n\left(\frac12C_{\phis}\, e_2, z\right) - \frac12C_{\phis}\, e_2(z) \right],
$$
thus passing to the limit, and recalling that $\ell''(z)=0$:
$$
- {M \over \phis(1)}\ \miu\ -A\, [\widehat{\phis}]''(0)\, + \, \frac12\, C_{\phis}\, [\widehat{\phis}]''(0)\ =\ -A\, [\widehat{\phis}]''(0)\, - \, {M \over \phis(1)},
$$
which gives a contradiction. The proof for the case of $t_{-}$ works analogously. Concerning the second part of the thesis, since $t_+$ is convex we have for any $h>0$:
$$
t_+(x+h)-2t_+(x)+t_+(x-h)\ \mau 0, \quad a<x<b.
$$
with $x+h$ and $x-h$ belonging to $(a,b)$. From this inequality, we immediately obtain:
$$
f(x+h)-2f(x)+f(x-h)\ \mau\ - C_{\phis}\, h^2.
$$
Similarly, using the convexity of $t_{-}$ we also find:
$$
f(x+h)-2f(x)+f(x-h)\ \miu\ C_{\phis}\, h^2.
$$
from which we deduce that $f \in Lip(2)$, with Lipschitz constant $C_{\phis}$.
\end{proof}
To conclude this part, we have to recall the well-known equivalence between the Lipschitz class $Lip(2)$ and the second-order Sobolev space $W^2_\infty(A)$, with $A=(a,b)$, which is crucial in order to establish a characterization of the saturation (Favard) classes of the NN operators $F_n$. Indeed, 
if we assume that $f \in W^2_\infty(A)$, where we recall that (see, \cite{DELO1}, p. 311):
$$
W^2_\infty(A)\, :=\, \left\{ f:A\to \R:\, f,f' \in AC(A),\, \mbox{with } |f''(x)| \miu M,\, \mbox{a.e.}\, x \in A  \right\}\footnote{$AC(A)$ denotes the usual set of all absolutely continuous functions.},
$$
since the second order Taylor formula with integral remainder holds (see p. 37 of \cite{DELO1}), and since $|f''(x)|\miu M$, a.e., in $[a,b]$, we can repeat the proof of Theorem \ref{voronovskaja-NN} (for $d=1$) on any compact $I_\delta=[a+\delta,b-\delta]$, $\delta>0$, in fact obtaining the inequality (\ref{ordine-locale-su-aperto}) with the constant $M>0$. This remark shows that also a sort of a converse implication of Theorem \ref{saturation-order-d=1} holds. 
\begin{example} \rm
Let $\sigma_\ell(x):=(1 + e^{-x})^{-1}$, $x \in \R$, be the well-known logistic function, which satisfies all the assumptions $(\Sigma i)$, $i=1,2,3$ (\cite{CCNP1,COPI2,CNP1}); in particular, $(\Sigma 3)$ is satisfied for every $\alpha>0$ since $\sigma_\ell$ decays exponentially as $x \to -\infty$. The Fourier transform of the corresponding $\phi_{\sigma_\ell}$ (which is, in fact, a cosine Fourier transform since $\phi_{\sigma_\ell}$ is even) is the following:
\be \label{eqn-39-new}
\widehat{\phi_{\sigma_\ell}}(v)\ :=\  {\pi\, \sin v \over \sinh(\pi v)},
\ee
with continuous extension at $v=0$.
This implies that $\widehat{\phi_{\sigma_\ell}}(2 \pi k) = 0$, for every $k \in \Z \setminus \left\{0 \right\}$, $\widehat{\phi_{\sigma_\ell}}(0) = 1$.  A proof of this fact has been given also in Theorem 1 of \cite{CAOCHEN2009}, exploiting integration methods of complex analysis. However, (\ref{eqn-39-new}) shows that $\sigma_\ell$ does not satisfy condition (\ref{ipotesi-trasformata}). Similar considerations can be done for the sigmoidal function $\sigma_h$ (see \cite{COSP3,CAOCHEN2015,KA1}) activated by the hyperbolic tangent function.

In order to give an example of sigmoidal function satisfying condition (\ref{ipotesi-trasformata}) we can proceed as follows.

We denote by:
$$
J_N(x)\ :=\ c_N\, \mbox{sinc}^{2 N}\left(\frac{x}{2 N\pi\gamma}\right), \hskip1cm x \in  \R,
$$
with $N \in \N$, $\gamma \geq 1$, and $c_N$ is a non-zero normalization coefficient, given by:
$$
c_N\ :=\ \left[ \int_{\R} \mbox{sinc}^{2 N}\left(\frac{u}{2 N \pi \gamma} \right) \, du \right]^{-1},
$$
the so-called {\em Jackson-type kernels} (see, e.g., \cite{BUNE}).
Here, the sinc-function $sinc(x)$ is defined as $\sin(\pi x)/ \pi x$, if $x \neq 0$, and $1$ if $x=0$. We now define the following sigmoidal function:
\be
\sigma_{J_N}(x)\, :=\, \int_{-\infty}^x J_N(t)\, dt, \quad x \in \R.
\ee
It is not difficult to see that $\sigma_{J_N}$ satisfies conditions $(\Sigma 1)-(\Sigma 3)$; in particular $\alpha$ is related to the parameter $N$. Thus, the density function corresponding to $\sigma_{J_N}$ can be written as follows:
$$
\phi_{\sigma_{J_N}}(x)\, =\, \frac12\, \int_{x-1}^{x+1} J_N(t)\, dt\, =\ \frac12\, \left( \chi_{[-1,1]}*J_N \right)(x),
$$
where the last term denotes the usual convolution product and $\chi_{[-1,1]}$ is the characteristic function of the interval $[-1,1]$. Thus, by the usual convolution theorem of Fourier Analysis, and recalling that
$supp\, \widehat{J_N}\, \subset [-1/\gamma, 1/\gamma]$ we immediately deduce that (\ref{ipotesi-trasformata}) holds, and all the previous results can be applied to $\sigma_{J_N}$.
\end{example}
%


\section{Saturation theorems for the Kantorovich NN by the solution of a hyperbolic PDE} \label{sec6}

For the Kantorovich NN operator, a local saturation theorem can be established thanks to the application of the generalized parabola techniques. For the sake of simplicity, we give the proof of the result below only in the bi-dimensional case $d=2$ and for the square $\widetilde{Q}:=[a,b]^2=[a,b]\times[a,b]$.
\begin{theorem} \label{th-saturation-local-KNN}
Let $\alpha>2$ {and assume that condition (\ref{ipotesi-trasformata}) holds for $s=1$}. Moreover, let $f \in C({\widetilde{Q}}) \cap C^1({\widetilde{Q}}^0)$ be fixed. If the function $f$ satisfies the following compatibility conditions:
\be \label{compatibility-condition-1}
f(x_1,a+\delta)=f(b-\delta,b-x_1+a), \quad x_1 \in [a+\delta, b-\delta],
\ee
\be \label{compatibility-condition-2}
f(a+\delta,x_2)=f(b-x_2+a,b-\delta), \quad x_2 \in [a+\delta, b-\delta],
\ee
{for a fixed $\delta>0$} and
\be \label{ipotesi-o-small-locale-KNN}
\sup_{{\tt x} \in {\widetilde{Q}}_{\delta}} |K_n(f, {\tt x} )-f({\tt x})|\ =\ o(n^{-1}), \quad as \quad n \to +\infty,
\ee
where ${\widetilde{Q}}_{\delta}=[a+\delta,b-\delta]^2$, then $f$ is a classical solution of the (hyperbolic) scalar divergence (or transport) equation on ${\widetilde{Q}}_\delta$:
\be \label{tesi-saturation-hyperbolic}
\sum_{i=1}^2 {\partial f \over \partial x_i}({\tt x})\, =\, 0, \quad {\tt x}=(x_1,x_2) \in {\widetilde{Q}}_\delta.
\ee
Vice-versa, {if $f \in C({\widetilde{Q}}) \cap C^1({\widetilde{Q}}^0)$ satisfies (\ref{tesi-saturation-hyperbolic}) then (\ref{compatibility-condition-1}), (\ref{compatibility-condition-2}) and (\ref{ipotesi-o-small-locale-KNN}) hold}.
\end{theorem}
\begin{proof}
Let $f \in C({\widetilde{Q}})\cap C^1({\widetilde{Q}}^0)$ satisfying (\ref{compatibility-condition-1}), (\ref{compatibility-condition-2}), and (\ref{ipotesi-o-small-locale-KNN}) for $\delta>0$, and consider the following boundary value problem:
\be \label{hyperbolic-boundary-value-problem}
\begin{cases}
\disp \sum_{i=1}^2 {\partial g \over \partial x_i}({\tt x})\, =\, 0, \quad {\tt x} \in {\widetilde{Q}}^0_\delta \\
g({\tt x})\, =\, f({\tt x}), \quad \quad {\tt x} \in \partial {\widetilde{Q}}_\delta.
\end{cases}
\ee
We know that in general the above problem might not have a solution. However, if we consider the initial value problem:
$$
\begin{cases}
\disp \sum_{i=1}^2 {\partial g \over \partial x_i}({\tt x})\, =\, 0, \quad \quad \quad \hskip1.1cm {\tt x} \in {\widetilde{Q}}^0_\delta \\
g(a+\delta,x_2)\, =\, f(a+\delta,x_2), \quad x_2  \mau a+\delta,
\end{cases}
$$
by the method of characteristics (see, e.g., Section 4 of \cite{salsa} and \cite{FJ1}), it is possible to see that it admits a unique classical solution $g$ (of class $C^1$, see also \cite{DIP1,BOGA}). However, this function $g$, thanks to the compatibility conditions (\ref{compatibility-condition-1}) and (\ref{compatibility-condition-2}) solves also the boundary value problem (\ref{hyperbolic-boundary-value-problem}).

Set now $G({\tt x})=g({\tt x}) - f({\tt x})$, ${\tt x} \in {\widetilde{Q}}_\delta$. It turns out that $G \in C({\widetilde{Q}}_\delta)$, and $G({\tt x})=0$ on $\partial {\widetilde{Q}}_\delta$. Reasoning as in the proof of Theorem \ref{th-saturation-NN-local}, $G$ must have a point of maximum (or minimum) in ${\widetilde{Q}}^0_\delta$ otherwise it must be $g \equiv f$ on ${\widetilde{Q}}_\delta$. Supposing that $G$ assumes its maximum of ${\tt x}_0 \in {\widetilde{Q}}^0_\delta$ (we can proceed similarly in the case of a minimum), by Theorem \ref{th-Ditzian-on-I} there exists $Q({\tt x}):= A\sum_{i=1}^2(x_i-a-\delta)^2+C$, with $A<0$, $C \in \R$, $Q({\tt x})>0$, $Q({\tt x})\mau G({\tt x})$ on ${\widetilde{Q}}_\delta$, $Q({\tt z})=G({\tt z})$, for a suitable ${\tt z} \in {\widetilde{Q}}_\delta^0$.

Let now consider the set ${\widetilde{Q}}_{\delta+\gamma} \subset {\widetilde{Q}}_\delta$, $\gamma>0$, and such that ${\tt z} \in {\widetilde{Q}}_{\delta+\gamma}$. Denoting by $\bar K_n$ the Kantorovich NN operators considered on the set ${\widetilde{Q}}_\delta$, as made in the proof of Theorem \ref{th-saturation-NN-local} for the operators $F_n$, it is possible to show (using (\ref{ipotesi-o-small-locale-KNN}) and the local Voronovskaja formula of Theorem \ref{voronovskaja-KNN}), that:
\be \label{auxiliary-order}
\sup_{{\tt x} \in {\widetilde{Q}}_{\delta+\gamma}} |\bar K_n(G, {\tt x} )-G({\tt x})|\ =\ o(n^{-1}), \quad as \quad n \to +\infty.
\ee
Now, we get:
$$
n\, \left[  \bar K_n(G, {\tt z} )-G({\tt z})  \right]\ \miu\ n\, \left[  \bar K_n(Q, {\tt z} )-Q({\tt z})  \right]\ =\ A\, \sum_{i=1}^2 n\, \left[  \bar K_n(\varphi_i, {\tt z} )-\varphi_i({\tt z})  \right],
$$
where $\varphi_i({\tt x})=(x_i-a-\delta)^2$, ${\tt x} \in {\widetilde{Q}}_\delta$, $i=1,2$. Passing to the limit in the above inequality, using (\ref{auxiliary-order}) and Theorem  \ref{voronovskaja-KNN} again, we obtain:
$$
0\ \miu\ A\,\{ (z_1-a-\delta) + (z_2-a-\delta)  \} \ \miu\ 2\, A\, \gamma\ <\ 0,
$$
hence a contradiction. This completes the proof. The second part of the thesis can be easily proved since, in general, all the solutions of (\ref{tesi-saturation-hyperbolic}) satisfy (\ref{compatibility-condition-1}) and (\ref{compatibility-condition-2}), while (\ref{ipotesi-o-small-locale-KNN}) follows by Theorem \ref{voronovskaja-KNN}. 
\end{proof}
\begin{remark} \rm
Note that, the compatibility conditions (\ref{compatibility-condition-1}) and (\ref{compatibility-condition-2}) are not surprising; indeed also in \cite{Zhou1} (see eqn. (1.4)) similar ones have been proved for the classical multivariate Kantorovich polynomials.
\end{remark}
\begin{remark} \rm
As for the case of the operators $F_n$, also for the Kantorovich NN operators we can get more accurate local approximation than the global ones. Moreover, we can also note that in the last part of the proof of Theorem \ref{th-saturation-local-KNN}, if we do not use the specific definition of the auxiliary function $Q$ explicitly pointed out in Lemma \ref{th-Ditzian-on-I}, we are not able to reach the desired result.
\end{remark}


\section{Semi-analytical inversion of the NN operators: a retrieval algorithm} \label{sec8}

In real-world applications, it may occur that when certain data are detected by measurements (by suitable sensors or measuring instruments), one can be interested not only in the resulting data, but also on the variables influencing these results. This interest arises primarily because some functional dependencies are known but can not be explicitly determined, leading to uncertainty about how variations in specific variables affect the observed outcomes.

Based on the above application problem, it could be useful to have a procedure available that, starting from the measured data allows one to estimate the value of certain target variables; clearly, this is a type of {\em inverse problem} that, in many applications, is referred to as a {\em retrieval problem}.

Retrieval algorithms are widely used in many engineering problems; for instance, in microwave remote sensing (RS) observations a wealth of methods have been proposed for the retrieval of essential climate variables (ECVs), see, e.g., \cite{OH1992}. Semi-empirical approaches are used to implement these kinds of algorithms; they are based on physical assumptions but rely on experimental trends, have limited portability across different measurement conditions and are particularly sensitive to noise. Moreover, widespread retrieval algorithms are based on data-driven techniques, which estimate the involved parameters directly from the experimental data via general statistical means, however, they are often computationally costly and require large amounts of training data. A method combining the simplicity and computational efficiency of semi-empirical approaches with the flexibility of data-driven techniques, while avoiding the use of extensive training datasets, is therefore desirable both as an alternative and as a complementary tool to existing retrieval schemes. 

Modeling by NN operators provides (approximate) functional expressions linking certain variables to suitable (also experimental) scalar values, even when the explicit functional dependencies are unknown. Hence, the theoretical inversion of the NN operator models can be useful to face the described retrieval problems.
\vskip0.2cm

Suppose now that we are interested in the estimation of a certain quantity $y \in \R$ which is dependent on $d$ physical variables, named $x_1, ..., x_d$, by an unknown functional relation of the form $y=f(x_1,...,x_d)$. 

Furthermore, we suppose that for a discrete uniformly spaced grid of values $\Sigma \subset \R^d$ , the corresponding values of $y$ are known. Denote these values by $y_{{\tt k}/n}$ to highlight that they have been measured correspondingly to the point ${\tt k}/n$ of the grid $\Sigma$; obviously, these measurements are typically affected by different sources of noise. In this context, one can propose the following approximate model describing the functional dependencies between the variable $x_i$, $i=1,...,d$, and $y$:
$$
(F_n f)({\tt x})\ =\ y,
$$  
for vectors ${\tt x}$ belonging to a given set containing the discrete grid $\Sigma$. 

  We are now interested in determining the contribution given by the $i$-th variable of a vector ${\tt x}$ in the determination of the value $y$; obviously we consider a vector that in general does not belong to the grid $\Sigma$, and we suppose that $y$ is known while the $i$-th variable of a vector ${\tt x}$ is unknown.

Moreover, recalling what has been observed in the Section 4 of \cite{CACO1}, if we assume that ${\tt x}$ is an inner point of $\Sigma$ (in the sense that it does not belong to the boundary points of the considered grid) the definition of the approximation operators $F_n$ can be simplified by omitting the denominator. More precisely, we get: 
$$
\sum_{{\tt k} \in {\cal J}_n} {y_{{\tt k}/n}}\, \Psi_\sigma ( n{\tt x} - {\tt k})\ =\ y,
$$
from which we can obtain:
\be \label{general-formulation}
\sum_{k_i \in {\cal I}^{[i]}_n}\phis(nx_i-k_i)\, \left[\sum_{{\tt k} \in {\cal J}^{[i]}_n} {y_{{\tt k}/n}}\, \Psi^{[i]}_\sigma ( n{\tt x}_{[i]} - {\tt k}_{[i]})\right]\ =\ y.
\ee
The equality in (\ref{general-formulation}) can be viewed as a nonlinear equation of unknown $x_i$ that can be solved by suitable numerical methods (see, e.g., \cite{KE2018}).

From the theory developed in \cite{COSP3}, it is also known that all the basic convergence results for the NN and the Kantorovich NN operators are valid even if assumption $(\Sigma 2)$ is not satisfied. In this case some alternative conditions have to be assumed, such as that the corresponding $\phis$ is a centered bell-shaped function.

  Based on the latter remark, we can here propose the following semi-analytical (retrieval) procedure valid in the special case of the so-called {\em ramp function},
which represents a (non differentiable) sigmoidal function for which the mentioned approximation results hold.

  We recall that the (sigmoidal) ramp function is defined as follows:
\begin{equation*}
\sigma_{R}(x)=\begin{cases}
0, & \text{ if } x<-\frac{1}{2}, \\
x+\frac{1}{2}, & \text{ if } -\frac{1}{2}\leq x\leq \frac{1}{2}, \\
1, & \text{ if } x>\frac{1}{2},
\end{cases}
\end{equation*} 
and the corresponding density function $\phi_{\sigma_{R}}$ can be explicitly written as follows:
$$
\phi_{\sigma_{R}}(x)\, =\, \begin{cases}
0, \hskip1.5cm |x|>{3/2},\\
1/2, \hskip1.1cm |x|<1/2,\\
\frac12 x +\frac34, \hskip0.6cm -3/2\miu x \miu-1/2,\\
-\frac12 x +\frac34 \hskip0.4cm 1/2\miu x\miu 3/2.
\end{cases}
$$
It is clear that $\phi_{\sigma_{R}}$ has compact support contained in $[-3/2,3/2]$ ({the length of its support is equal to} 3).
Choosing $\sigma=\sigma_R$ in (\ref{general-formulation}) and recalling the specific form of  the support of $\phi_{\sigma_{R}}$ we can obtain:
$$
\sum_{k_i = \bar k -1}^{\bar k +1}\phi_{\sigma_R}(nx_i-k_i)\, \left[\sum_{{\tt k} \in {\cal J}^{[i]}_n} {y_{{\tt k}/n}}\, \Psi^{[i]}_{\sigma_R} ( n{\tt x}_{[i]} - {\tt k}_{[i]})\right]\ =\ y,
$$
for a certain $\bar k$ which is strictly greater than $\lceil  n a_i\rceil$ and strictly less than $\lfloor nb_i \rfloor$. Since the above sum with respect to $k_i$ reduces only to three terms, they must be distributed as:
\be \label{distribution-nodes}
nx_i-\bar k-1 \in \left[-\frac32, \frac12\right), \quad nx_i-\bar k \in \left[-\frac12, \frac12\right], \quad nx_i-\bar k+1 \in \left(\frac12, \frac32\right].
\ee
Hence we get:
$$
\phi_{\sigma_R}(nx_i-\bar k -1)\, \left[\sum_{{\tt k} \in {\cal J}^{[i]}_n} {y}_{{\tt k}^{[i]}_{[\bar k +1]} / n}\, \Psi^{[i]}_{\sigma_R} ( n{\tt x}_{[i]} - {\tt k}_{[i]})\right]\ 
$$
$$
+\ \phi_{\sigma_R}(nx_i-\bar k)\, \left[\sum_{{\tt k} \in {\cal J}^{[i]}_n} {y}_{{\tt k}^{[i]}_{[\bar k]} / n}\, \Psi^{[i]}_{\sigma_R} ( n{\tt x}_{[i]} - {\tt k}_{[i]})\right]\ 
$$
$$
+\ \phi_{\sigma_R}(nx_i-\bar k +1)\, \left[\sum_{{\tt k} \in {\cal J}^{[i]}_n} {y}_{{\tt k}^{[i]}_{[\bar k -1]} / n}\, \Psi^{[i]}_{\sigma_R} ( n{\tt x}_{[i]} - {\tt k}_{[i]})\right]\ =\ y
$$
where by the symbol ${\tt k}^{[i]}_{[a]}$ we denote the $d$-dimensional vector having in the $i$-th position the specific value $a \in \R$. Now, using (\ref{distribution-nodes}) and the definition of $\phi_{\sigma_R}$ we obtain:
$$
\left[\frac12(nx_i-\bar k -1)+\frac34\right]\, \left[\sum_{{\tt k} \in {\cal J}^{[i]}_n} {y}_{{\tt k}^{[i]}_{[\bar k +1]} / n}\, \Psi^{[i]}_{\sigma_R} ( n{\tt x}_{[i]} - {\tt k}_{[i]})\right]\ 
$$
$$
+\ \frac12 \left[\sum_{{\tt k} \in {\cal J}^{[i]}_n} {y}_{{\tt k}^{[i]}_{[\bar k]} / n}\, \Psi^{[i]}_{\sigma_R} ( n{\tt x}_{[i]} - {\tt k}_{[i]})\right]\ 
$$
$$
+\ \left[ -\frac12 (nx_i-\bar k +1) + \frac34  \right]\, \left[\sum_{{\tt k} \in {\cal J}^{[i]}_n} {y}_{{\tt k}^{[i]}_{[\bar k -1]} / n}\, \Psi^{[i]}_{\sigma_R} ( n{\tt x}_{[i]} - {\tt k}_{[i]})\right]\ =\ y.
$$
From the above linear equation we can easily extract the retrieved variable $x_i$, namely: 
$$
x_i\, =\ n^{-1}\, \left\{  \sum_{{\tt k} \in {\cal J}^{[i]}_n}{y}_{{\tt k}^{[i]}_{[\bar k +1]} / n}\, \Psi^{[i]}_{\sigma_R} ( n{\tt x}_{[i]} - {\tt k}_{[i]})\, -\,  \sum_{{\tt k} \in {\cal J}^{[i]}_n} {y}_{{\tt k}^{[i]}_{[\bar k -1]} / n}\, \Psi^{[i]}_{\sigma_R} ( n{\tt x}_{[i]} - {\tt k}_{[i]})   \right\}^{-1} \times
$$
\be \label{eqn50}
\left\{ 2\, y\, +\, \left(  \bar k - \frac12  \right)\, \left[   \sum_{{\tt k} \in {\cal J}^{[i]}_n} {y}_{{\tt k}^{[i]}_{[\bar k +1]} / n}\, \Psi^{[i]}_{\sigma_R} ( n{\tt x}_{[i]} - {\tt k}_{[i]}) \right] \right.
\ee
$$
-\, \left[\sum_{{\tt k} \in {\cal J}^{[i]}_n} {y}_{{\tt k}^{[i]}_{[\bar k]} / n}\, \Psi^{[i]}_{\sigma_R} ( n{\tt x}_{[i]} - {\tt k}_{[i]})\right]
\left. -\, \left(  \bar k + \frac12  \right)\, \left[\sum_{{\tt k} \in {\cal J}^{[i]}_n} {y}_{{\tt k}^{[i]}_{[\bar k -1]} / n}\, \Psi^{[i]}_{\sigma_R} ( n{\tt x}_{[i]} - {\tt k}_{[i]})\right]  \right\},
$$
assuming that the following term is non-zero:
\be \label{multivariate-retrieval-formula}
\sum_{{\tt k} \in {\cal J}^{[i]}_n} {y}_{{\tt k}^{[i]}_{[\bar k +1]} / n}\, \Psi^{[i]}_{\sigma_R} ( n{\tt x}_{[i]} - {\tt k}_{[i]})\, -\,  \sum_{{\tt k} \in {\cal J}^{[i]}_n} {y}_{{\tt k}^{[i]}_{[\bar k -1]} / n}\, \Psi^{[i]}_{\sigma_R} ( n{\tt x}_{[i]} - {\tt k}_{[i]}).
\ee
The retrieval formula given in (\ref{eqn50}) can be easily implemented using suitable programming languages, such as MATLAB.
\vskip0.2cm

  For the sake of completeness, we recall that the present case is strictly related to the so-called \textit{Rectified Linear Unit} (ReLU) activation function (\cite{KOLOM,LABATE,DENO}) defined by $\psi_{\text{ReLU}}(x)=(x)_{+}:= \max \{x, 0\}$, $x\in\R$ (\cite{RUT}). More precisely, $\phi_{\sigma_R}$ can be expressed in terms of the ReLU as follows:
\begin{equation*}
\phi_{\sigma_R}(x)=\frac{1}{2}\left[\psi_{\text{ReLU}}\left(x+\frac{3}{2}\right)-\psi_{\text{ReLU}}\left(x+\frac{1}{2}\right)-\psi_{\text{ReLU}}\left(x-\frac{1}{2}\right)+\psi_{\text{ReLU}}\left(x-\frac{3}{2}\right)\right],
\end{equation*}
$x\in\R$. {Note that, if we consider the activation functions $\sigma_{M_n}$ generated by the well-known central B-spline of order $n$ (for the definition see \cite{CO-reluk}) the corresponding $\phi_{\sigma_{M_n}}$ can be viewed as a combination of $ReLU^k$ activation functions (see, e.g., \cite{YAN24,YAN25,SHE}), which are other typical examples of activation functions of recent NN models.}
\begin{example} \rm
Here we briefly give two numerical experiments (both for functions of one and two variables) of the above retrieval algorithm. Let us consider the function $f_1(x):=x^2+1$, on the interval $[0,1]$, {that we used to compute the values $y_{k/n}$ in (\ref{eqn50})}. In Table \ref{test1} we have the retrieved variable $x$ by means of the considered NN operators $F_n$ for $n=5, 10, 15$.
\begin{table}[h!]
\centering
\begin{tabular}{ccccc}
\hline
${\footnotesize \bf n}$ & \vline &  {\footnotesize \bf retrieved nodes} &   & {\footnotesize \bf max-error} \\
\hline
 {\footnotesize \bf 5} \!\!  & \vline  & {\scriptsize 0.1500  \,  0.3750 \,   0.5833 \,  0.7875  } & \vline & {\footnotesize 0.0500 } \\ 
{\footnotesize \bf 10} \!\! & \vline  & {\scriptsize  0.0750  \,  0.1875 \,   0.2917  \,  0.3938  \,  0.4950  \,  0.5958  \,  0.6964 \,   0.7969  \,  0.8972 } & \vline & {\footnotesize 0.0250}  \\
{\footnotesize \bf 15} \!\! & \vline  &  {\scriptsize  0.0500    0.1250    0.1944    0.2625    0.3300    0.3972    0.4643    0.5313    0.5981   
 0.6650 0.7318    0.7986    0.8654    0.9321} \!\!  & \vline & {\footnotesize 0.0167}  \\
\hline
\end{tabular}
\caption{The retrieved values for the variable $x$ from the application of the semi-analytical formula (\ref{eqn50}) for $d=1$ to $f_1$.}
\label{test1}
\end{table}
The max-errors in the last column in Table \ref{test1} have been computed using the analytical nodes as reference values.

As a second example we can consider the bivariate function $f_2(x_1,x_2):=x_1^3+x_2+1$, on the square $\widetilde{Q}=[0,1]^2$, {that we used to compute the values $y_{{\tt k}/n}$ in (\ref{eqn50})}. In Table \ref{test2} we reported the retrieved grid of nodes by means of the considered NN operators $F_n$ for $n=5, 10$.
\begin{table}[h!]
\centering
\begin{tabular}{ccccc}
\hline
 & \vline &  {\footnotesize \bf retrieved nodes for each variable} &   & {\footnotesize \bf max-error } \\
\hline
${\scriptsize \bf n = 5}$ & \vline &   &   &  \\
\hline
 {\footnotesize $x_1$} \!\!  & \vline  & {\scriptsize 0.1290    0.3794    0.6205    0.8563 } & \vline & {\footnotesize 0.0710 } \\ 
{\footnotesize $x_2$} \!\! & \vline  & {\scriptsize  0.2149    0.4434    0.6720    0.9006 } & \vline & {\footnotesize 0.1006}   \\
\hline
${\scriptsize \bf n = 10}$ & \vline &    &   &   \\
\hline
 {\footnotesize $x_1$} \!\!  & \vline  & {\scriptsize 0.0645    0.1897    0.3102    0.4282    0.5448    0.6607    0.7762    0.8913    1.0063 } & \vline & {\footnotesize 0.1063 } \\ 
{\footnotesize $x_2$} \!\! & \vline  & {\scriptsize  0.1126    0.2269    0.3411    0.4554    0.5697    0.6840    0.7983    0.9126    1.0269 } & \vline & {\footnotesize 0.1269}   \\
\hline
\end{tabular}
\caption{The retrieved grid of nodes from the application of the semi-analytical formula (\ref{eqn50}) for $d=2$ to $f_2$.}
\label{test2}
\end{table}

  It is clear that, from the numerical results of Table \ref{test2} one can construct the full grid of nodes by computing the simple Cartesian product between the nodes of the variables $x_1$ and $x_2$. 
\end{example}


\section{Final remarks and open problems}

In the present paper we have thoroughly investigated the saturation and inverse problems for two families of multivariate NN-type operators. The proposed techniques of proof are mainly constructive, and some of the established results are linked to the solution of elliptic and hyperbolic PDEs. Although the study presented here is comprehensive, we still have some open problems. For instance, concerning the operators $K_n$ the situation is really delicate. It would be desirable to remove from the saturation theorem of Section \ref{sec6} the assumed compatibility conditions. Global/Local inverse theorems with respect to the $\sup$-norm are missing, as are analogous results in the $L^p$-norm setting. As observed by Zhou (\cite{Zhou1}), multivariate inverse problems for Kantorovich type operators require specific techniques (also related to the theory of distribution or to interpolation theory) that we plan to explore in future work. As fundamental results for Kantorovich type operators one can see \cite{TOT1,TOT2}.

Finally, we emphasize that the semi-analytical problem discussed in Section \ref{sec8} represents only a preliminary attempt to tackle a topic that could be fundamental in the development of the research project RETINA (mentioned below in the funding section of the paper) in which theoretical retrieval methods for the application to remote sensing data should be applied for extracting information on the health status of land surfaces, and to help the supporting forecasts of desertification and/or flood risk.   


\section*{Acknowledgments}

{\small The author would like to thank the anonymous referees for their careful reading of the manuscript and for their valuable comments and suggestions, which have contributed to improving the paper. 

  The author is member of the Gruppo Nazionale per l'Analisi Matematica, la Probabilit\`a e le loro Applicazioni (GNAMPA) of the Istituto Nazionale di Alta Matematica (INdAM), of the network RITA (Research ITalian network on Approximation), and of the UMI (Unione Matematica Italiana) group T.A.A. (Teoria dell'Approssimazione e Applicazioni). 
}

\section*{Funding}

{\small The author has been supported within the PRIN 2022 PNRR: ``RETINA: REmote sensing daTa INversion with multivariate functional modeling for essential climAte variables characterization", funded by the European Union under the Italian National Recovery and Resilience Plan (NRRP) of NextGenerationEU, under the Italian Ministry of University and Research (Project Code: P20229SH29, CUP: J53D23015950001).
}
\section*{Conflict of interest/Competing interests}

{\small The author declares that he has no conflicts of interest or competing interests.}

\section*{Availability of data and material and Code availability}

{ \small Not applicable.}

%

\end{document}